\documentclass[a4paper]{article}
\usepackage{amssymb,amsmath,amsthm}
\numberwithin{equation}{section}

\newcommand{\Q}{\mathbb{Q}}

\newcommand{\R}{\mathbb{R}}

\newcommand{\N}{\mathbb{N}}
\newcommand{\Z}{\mathbb{Z}}

\newtheorem*{thm}{Theorem}
\newtheorem{lem}{Lemma}
\newtheorem{pro}{Proposition}

\theoremstyle{definition}
\newtheorem*{ack}{Acknowledgements}

\renewcommand{\mod}[1]{\hspace{-2.9mm}\pmod{#1}}
\newcommand{\x}{{\bf x}}
\newcommand{\e}{\emph}
\newcommand{\rom}{\mathrm}
\newcommand{\bfP}{\mathbb{P}}
\newcommand{\A}{\mathbb{A}}

\newcommand{\ma}{\mathbf}
\newcommand{\ben}{\begin{enumerate}}
\newcommand{\een}{\end{enumerate}}
\newcommand{\eit}{\begin{itemize}}
\newcommand{\beq}{\begin{equation}}
\newcommand{\eeq}{\end{equation}}

\newcommand{\ve}{\varepsilon}

\newcommand{\mcal}{\mathcal}
\newcommand{\lab}{\label}
\newcommand{\al}{\alpha}

\newcommand{\be}{\beta}

\newcommand{\sfl}{\mathsf{\Lambda}}

\newcommand{\hcf}{\mathrm{gcd}}
\newcommand{\colt}[2]{\genfrac{}{}{0pt}{1}{#1}{#2}}

\newcommand{\eqm}[2]{\equiv #1 \pmod{#2}}

\newcommand{\vt}{\vartheta}

\renewcommand{\leq}{\leqslant}
\renewcommand{\geq}{\geqslant}
\renewcommand{\d}{\mathrm{d}}

\begin{document}

\title{The density of rational points on a certain singular cubic surface}
\author{T.D. Browning\\
\small{\emph{School of Mathematics, Bristol University, Bristol BS8 1TW}}\\ 
\small{t.d.browning@bristol.ac.uk}}
\date{}

\maketitle

\begin{abstract}
We show that the number of non-trivial rational 
points of height at most $B$, which lie on the cubic surface
$x_1x_2x_3=x_4(x_1+x_2+x_3)^2$,  has order of magnitude $B(\log B)^6$.
This agrees with Manin's conjecture.
\end{abstract}

\section{Introduction}

The purpose of this paper is to investigate the distribution of
rational points on the singular cubic surface $X \subset \bfP^3$,
given by the equation
$$
x_1x_2x_3=x_4(x_1+x_2+x_3)^2.
$$
This surface has a unique singular point $[0,0,0,1]$ which is of type
$D_4$, and contains 
precisely $6$ lines  \cite[Lemma $4$]{b-w}.  These lines are all
defined over $\Q$ and are given by
$$
x_i=x_4=0, \quad 
x_i=x_j+x_k=0,
$$
for $\{i,j,k\}=\{1,2,3\}$.
We shall denote by $U \subset X$ the open subset formed by deleting the lines
from $X$.   

Given a rational point $x=[x_1,x_2,x_3,x_4] \in \bfP^3(\Q)$
such that $x_1,x_2,x_3,x_4$ are relatively prime integers, let
$H(x)=\max_{1\leq i \leq 4}|x_i|$ denote its 
anticanonical height, metrized by the choice of norm $\max_{1\leq i \leq 4}|x_i|$.
Then for any $B \geq 1$, we shall be concerned with
estimating the quantity
$$
N_{U,H}(B)= \#\{x \in U\cap\bfP^{3}(\Q): H(x) \leq B\}.
$$
Manin \cite{frankemanintschinkel89}  
has provided a very general conjecture concerning the
distribution of rational points on Fano varieties.  In our case it
predicts that there exists a positive constant $c_{X,H}$ such that
$$
N_{U,H}(B) \sim  c_{X,H} B (\log B)^6, 
$$
as $B \rightarrow \infty$.  Here the exponent of $\log B$ is one less
than the rank of the Picard group $\rom{Pic}\tilde{X}$ of $\tilde{X}$,
where $\tilde{X}$ denotes the minimal desingularisation of $X$.
In fact this sort of asymptotic formula is conjectured to hold
for any cubic surface with canonical singular locus. Although there
has been increasing interest
in Manin's conjecture for cubic surfaces, it has only been completely settled 
in particularly simple cases such as the toric variety 
$$
x_1x_2x_3=x_4^3.
$$ 
This can be found in the work of 
la Bret\`{e}che \cite{b}, for example.

More recently, Heath-Brown \cite{cayley} has 
established upper and lower bounds for the density of
non-trivial rational
points on the Cayley cubic surface
$$
\frac{1}{x_1}+\frac{1}{x_2}+\frac{1}{x_3}+\frac{1}{x_4}=0,
$$
which agree with Manin's conjecture.   
This is a cubic surface containing  four $A_1$ singularities, which is
the maximal number of singularities that a non-ruled cubic surface can
have.  The principal tool in Heath-Brown's work is a passage to the universal
torsor above the minimal desingularisation of the Cayley cubic.  Originally
introduced by Colliot-Th\'el\`ene and Sansuc to aid in
the study of the Hasse principle and weak approximation, universal torsors
were first used by Peyre \cite{peyre} and Salberger \cite{s} in the
context of counting rational 
points of bounded 
height.  
After establishing a bijection between the rational points on the
Cayley cubic and the integer points on the universal torsor, 
which in this setting is 
given explicitly by nine equations in thirteen variables,
Heath-Brown proceeds
by applying methods from the geometry of numbers to count 
integer solutions to certain ternary linear equations.

Our present work 
is largely inspired by Heath-Brown's treatment of the Cayley cubic surface.
We are now in a position to state our main result.

\begin{thm}
We have
\begin{equation}\lab{theorem}
B (\log B)^6 \ll N_{U,H}(B) \ll B (\log B)^6.
\end{equation}
\end{thm}

Of the two bounds in our theorem, the lower bound $N_{U,H}(B)\gg B
(\log B)^6$ is routine.  It will follow from relatively minor
adjustments to Heath-Brown's treatment of the 
Cayley cubic.
Establishing the upper bound in (\ref{theorem}), however, is by
far the most challenging component of this paper and the proof has two fundamental
ingredients.  The first is a translation of the problem to the
universal torsor above $\tilde{X}$, which in this setting has
the affine embedding  
\begin{equation}\lab{ut*} 
s_0s_1s_2s_3u_1u_2u_3 = y_1u_1s_1^2 + y_2u_2s_2^2 + y_3u_3s_3^2.
\end{equation}
This has been calculated by Hassett and Tschinkel
\cite[\S $4$]{hassetttschinkel03}, although we shall present our own
deduction of this equation in \S \ref{ut} below.
The universal torsor can be thought of as
serving  to encode factorisation information about the integer solutions
to the original 
equation.  In practical terms, it allows us to work with a 
larger number of variables, all of which are smaller in
modulus than the original variables $x_1,x_2,x_3,x_4$.  
The second main ingredient in our proof of the upper bound involves
studying the distribution of integer solutions to the equation
obtained by setting $s_0=0$ in (\ref{ut*}).  This is the focus of \S
\ref{core-x} and relies upon lattice methods to count 
integer solutions to ternary linear and quadratic forms.  
It seems worthwhile highlighting the fact that this need to consider
the contribution from quadratic equations marks a 
significant departure from Heath-Brown's treatment of the Cayley cubic.
In particular, we shall need to pay careful attention to the fact that
almost all ternary quadratic forms  don't have an integer solution.

Our work draws upon a diverse range of techniques.
In addition to
the geometry of numbers used to study linear and
quadratic forms in \S \ref{diag}, we make use of the 
large sieve inequality 
and real character sum estimates due to Heath-Brown  and
P\'olya--Vinogradov in \S \ref{sieve-in}.

\begin{ack}
This problem was posed by Yuri Tschinkel during the 
American Institute of Mathematics workshop ``Rational and integral
points on higher-dimensional varieties'' in December, 2002.  
The author is very grateful to Professor Heath-Brown and Professor
Tschinkel for several useful conversations relating to the subject of
this work.  Special thanks are due to the anonymous referee for his careful reading of the
manuscript and numerous pertinent remarks.
While working on this paper, the 
author was supported at Oxford University by EPSRC grant number GR/R93155/01.
\end{ack}

\section{Preliminary estimates}\lab{prelim}

We begin by introducing some conventions regarding our choice of
notation.  Throughout this paper 
the letters $i,j,k$ will denote generic distinct indices from the set $\{1,2,3\}$.
We shall use $\N$ to denote the set of positive integers, and for any
$n \geq 2$ it will be convenient to let $Z^n$ denote the 
set of primitive vectors in $\Z^n$,
where $\ma{v} \in \Z^n$ is
said to be primitive if 
$
\hcf(v_1,\ldots,v_n)=1.
$
Similarly, we let $N^n$ denote the set  of primitive vectors in $\N^n$.
Furthermore, we let  $Z_*^n$  denote the subset
of $\ma{v}\in Z^n$ for which $v_1\cdots v_n\neq 0$.  
Upon writing 
$$
F(\x)=x_1x_2x_3-x_4(x_1+x_2+x_3)^2, 
$$ 
it therefore follows that 
\beq\lab{n(b)}
N_{U,H}(B)= \frac{1}{2}\#\Big\{\x\in Z_*^4: \max_{1\leq i\leq 4}|x_i|\leq B, ~F(\x)=0\Big\},
\eeq
since $\x$ and $-\x$ represent the same point in $\bfP^3$.  
It will be convenient to collect together some technical results that will be
useful to us.

\subsection{The geometry of numbers and ternary forms}\lab{diag}

Several of our arguments will involve estimating the number 
of primitive integer solutions to certain  ternary homogeneous polynomial
equations, which lie in lopsided regions.
In the case of linear equations, such 
an estimate is provided by the following result of 
Heath-Brown  \cite[Lemma $3$]{h-b84}.

\begin{lem}\lab{line}
Let $\ma{h} \in Z^3$ and let
$W_i >0$.  Then the number of $\ma{w} \in Z^3$ for which 
$\sum_{i=1}^3h_i w_i=0,
$ 
and $|w_i| \leq W_i$,
is 
$$
\leq 4+12 \pi \frac{W_1W_2W_3}{\max |h_i|W_i}.
$$
\end{lem}

We shall also need a result which handles the corresponding problem
for diagonal quadratic equations.
For this we turn to the following 
result, in which 
$\omega(n)$ denotes the number of distinct prime
factors of $n \in \N$.  

\begin{lem} \lab{box}
Let $\ma{g, h} \in Z_*^3$, with $g_1g_2g_3$ square-free, and let
$W_i >0$.  Then the number of $\ma{w} \in Z^3$ for which 
$
\sum_{i=1}^3 g_i h_i w_i^2=0,
$ 
and $|w_i| \leq W_i$, is 
$$
\ll \left(1+ \sqrt{\frac{W_1W_2W_3 D_\ma{g,h}^{3/2}}{|h_1h_2h_3|}}\right)2^{\omega(h_1h_2h_3)},
$$
where $D_\ma{g,h}$ is the product of greatest common divisors
$$
\hcf(h_1h_2,h_1h_3,h_2h_3)\hcf(g_1,h_2h_3)\hcf(g_2,h_1h_3)\hcf(g_3,h_1h_2).
$$
\end{lem}

Lemma \ref{box} will follow from a rather straightforward modification to the proof
of Heath-Brown's \cite[Theorem 2]{h-b97}. 
In fact Heath-Brown establishes a version of Lemma \ref{box} with
$\ma{g}=(1,1,1)$ and 
$d_3(|h_1h_2h_3|)$ in place of $2^{\omega(h_1h_2h_3)}$, where 
$d_k(n)$ denotes the number of representations of $n$ as a product of
$k$ positive integers, for any $k,n \in \N$.  It is perhaps worth pointing out that whereas
$d_3(n)$ has average order $\frac{1}{2}(\log n)^2$, the function $2^{\omega(n)}$
has average order $\zeta(2)^{-1}\log n$.  This saving plays an important
role in our work.

In order to prove Lemma \ref{box} we  recall that 
the original idea behind the proof of \cite[Theorem 2]{h-b97}  
is to view the equation $\sum_{i=1}^3 g_ih_iw_i^2=0$ as a collection of
lattice conditions upon the  solutions $\ma{w} \in \Z^3$.
Let $p$ be any prime divisor of $h_1h_2h_3$, and assume without loss
of generality that 
$$
0 \leq \nu_p(h_1)\leq \nu_p(h_2) \leq \nu_p(h_3),
$$
where $\nu_p(n)$ denotes the $p$-adic order of any $n \in \N$.  
In particular it follows that $\nu_p(h_1)=0$, since $\ma{h}$ is primitive.
We shall only consider here the case in which $p$
is an odd prime. The case $p=2$ is handled along similar lines.
Since $g_1g_2g_3$ is square-free, we may write
$$
g_1=p^{\al_1}g_1', \quad g_2=p^{\al_2}g_2', 
\quad g_3=p^{\al_3}g_3',
$$
for $(\al_1,\al_2,\al_3) \in \{(0,0,0),(1,0,0),(0,1,0),(0,0,1)\}$ and
$p \nmid g_i'$.  Similarly, we write
$$
h_2=p^{\be_2}h_2', \quad h_3=p^{\be_3}h_3', 
$$
for $p \nmid h_2'h_3'$ and $\be_3\geq \be_2 \geq 1$.  
Then one proceeds by considering solutions to the equation
\beq\lab{cong1}
p^{\al_1}g_1'h_1w_1^2+p^{\al_2+\be_2}g_2'h_2'w_2^2 +p^{\al_3+\be_3}g_3'h_3'w_3^2=0.
\eeq

Suppose for the moment that we are examining solutions
$(u,v,w) \in \Z^3$ to the equation
\beq\lab{s-t}
au^2+p^\sigma bv^2 +p^{\tau}cw^2=0,
\eeq
for $0 \leq \sigma \leq \tau$ and $p \nmid abc$. Then arguing along
similar lines to the proof of \cite[Theorem 2]{h-b97}, we sketch how
this implies that $(u,v,w)$ lies on one of 
at most $2$ sublattices of $\Z^3$, each of
determinant $p^{\delta(\sigma,\tau)}$, where
\beq\lab{det_p}
\delta(\sigma,\tau) =
\left\{
\begin{array}{ll}
{(\sigma+\tau)-3\sigma/2}, &   \mbox{$\sigma$ even,}\\
{(\sigma+\tau)-[3\sigma/2]+1}, &   \mbox{$\sigma$ odd.}
\end{array}
\right.
\eeq
Suppose first that $\sigma=2s$ is even.  Then (\ref{s-t}) implies that
$p^s \mid u$.  By writing $u=p^su'$, and considering the corresponding
congruence $au'^2+bv^2 \eqm{0}{p^{\tau-\sigma}}$, we therefore deduce
that $(u,v,w)$ lies on one of at most $2$ integer lattices, each  of
determinant $p^{s+\tau-\sigma}=p^{\delta(\sigma,\tau)}$.  Suppose
now that $\sigma=2s+1$ is odd.
In view of (\ref{s-t}) we may again write $u=p^su'$, and consider the corresponding
congruence $au'^2+p b v^2 \eqm{0}{p^{\tau-2s}}$.  Since $\tau-2s \geq
1$,  
we may clearly write $u'=pu''$, and so consider solutions to the equation 
$$
p a ({u''})^2+b v^2 +p^{\tau-\sigma}cw^2=0.
$$    Now either
$\tau-\sigma=0$, or else we may write $v=pv'$ and consider the
equation $a ({u''})^2+pb v'^2 +p^{\tau-\sigma-1}cw^2=0$.  In the
former case we conclude that $(v,w)$ lies on one of at most two
integer sublattices of determinant $p$.  But then 
$(u,v,w)$ lies on one of at most $2$ integer lattices, each  of
determinant $p^{s+2}=p^{\delta(\sigma,\tau)}$. 
In the latter case we have $\tau-\sigma \geq
1$, and  we proceed inductively.  Thus either $\tau-\sigma=1$, in which case we
deduce that $(u,v,w)$ lies on one of $2$ integer lattices of
determinant $p^{s+3}=p^{\delta(\sigma,\tau)}$, or else
$\tau-\sigma\geq 2$ and we can repeat the process.  Since this process
clearly terminates we therefore deduce that whenever $\sigma$ is
odd, $(u,v,w)$ lies on one of at most $2$ integer lattices, each  of
determinant $p^{\delta(\sigma,\tau)}$.

Returning to (\ref{cong1}), our goal is to show that
$\ma{w}$ lies on one of at most $2$ integer sublattices of $\Z^3$,
each of determinant 
\beq\lab{test}
\geq p^{\be_2+\be_3-[3(\al_1+\al_2+\al_3+\be_2)/2]}.
\eeq
In view of the existing proof of \cite[Theorem 2]{h-b97}, this will suffice to
establish Lemma \ref{box} since the inequalities $\be_3 \geq \be_2
\geq 1$ imply that 
$$
\nu_p(h_1h_2h_2)=\be_2+\be_3, \quad \nu_p(D_\ma{g,h})=\al_1+\al_2+\al_3+\be_2.
$$ 
Suppose first that $(\al_1,\al_2)=(0,0)$.  
Then our work above shows that $\ma{w}$ lies on one of at most $2$
integer sublattices of $\Z^3$, each of determinant 
$$
p^{\delta(\be_2,\al_3+\be_3)} \geq p^{\delta(\be_2,\be_3)}.
$$
This is plainly satisfactory for (\ref{test}), by (\ref{det_p}).
Suppose now that $(\al_1,\al_2,\al_3)=(1,0,0)$.   
If $\be_2=0$, then it is not hard to conclude that
$\ma{w}$ lies on one of at most $2$ lattices of determinant
$p^{\be_3}$, which is also satisfactory.  
If now $\be_2 \geq 1$ we obtain an equation of the shape 
 (\ref{s-t}), with  $\sigma=\be_2-1$ and
$\tau=\be_3-1$.  Thus we obtain at most
$2$ integer lattices, each of determinant
$p^{\delta(\be_2-1,\be_3-1)}.$   It is easily checked that this
quantity is bounded below by (\ref{test}).
Finally we suppose that  $(\al_1,\al_2,\al_3)=(0,1,0)$.  
In this case we
again obtain an equation of the form (\ref{s-t}).  Suppose first that
$\be_2<\be_3$, so that we may take $\sigma=\be_2+1$ and
$\tau=\be_3$ in (\ref{s-t}).  But then it easily follows that
$\ma{w}$ lies on one of at most
$2$ integer lattices, each of determinant
$$
p^{\delta(\be_2+1,\be_3)}\geq p^{\be_2+\be_3-[3(1+\be_2)/2]}.
$$ 
Alternatively, if $\be_2=\be_3=\be$ say, then we take  
$\sigma=\be$ and $\tau=\be+1$ in (\ref{s-t}), thereby deducing that 
$\ma{w}$ lies on one of at most $2$ lattices, each of determinant
$$
p^{\delta(\be,\be+1)}\geq p^{2\be-[3(1+\be)/2]}.
$$ 
This completes the proof of Lemma \ref{box}.

\subsection{Solubility of quadratic forms}\lab{sieve-in}

In addition to considering
the density of integer solutions to 
diagonal quadratic equations, as in the previous section, we shall  also
need to consider how often such an equation has at least one
non-trivial integer solution.   
Let $\ma{a}\in Z_*^3$, let $Y_1,Y_2,Y_3 \geq 1$, and let $H \in \N$.
We shall write $T(\ma{Y};\ma{a},H)$ to denote
the set of pairwise coprime non-zero integers $y_1,y_2,y_3$ such that $|y_i|\leq Y_i$,
and 
\beq\lab{divH}
\hcf(a_iy_i,a_jy_j)\mid H, 
\eeq
and such that the
equation
$$
a_1y_1x_1^2+a_2y_2x_2^2+a_3y_3x_3^2=0
$$
has a non-zero solution $\x \in \Z^3$ with $\hcf(x_i,x_j)=1$.  
If $\ma{y} \in T(\ma{Y};\ma{a},H)$ then it necessarily follows that 
\beq\lab{jac-1}
\Big(\frac{-a_ia_jy_iy_j}{p}\Big)\neq -1 \quad \mbox{for all odd
  primes $p \mid a_ky_k$},
\eeq
where $(\frac{n}{p})$ denotes the Legendre symbol for any $n \in \Z$
and odd prime $p,$ and as usual $\{i,j,k\}$ denotes any permutation of the set $\{1,2,3\}$.
Define the arithmetic function
\beq\lab{phi}
\vt(n)=\prod_{p \mid n} \Big(1+\frac{1}{p}\Big),
\eeq
for any $n \in \N$.   
We shall proceed under the assumption that 
\beq\lab{scat}
Y_1 \leq Y_2 \leq Y_3.
\eeq
With this in mind  the goal of this section is to establish the following estimate.

\begin{pro}\lab{guo}
Let $\ve>0$.  Then we have
$$
\sum_{\ma{y}\in T(\ma{Y};\ma{a},H)} 2^{\omega(y_1y_2y_3)}\ll_\ve
\vt(a_1a_2)H^\ve \Big(Y_1Y_2Y_3
+(Y_1Y_2)^{1/2+\ve}Y_3m_\ve(\ma{a,Y})
\Big),
$$
with
\beq\lab{m}
m_\ve(\ma{a,Y})=\min\{|a_1a_2|,Y_3\}^{\ve}+\log Y_3,
\eeq
\end{pro}

Before proceeding with the proof of Proposition \ref{guo}, it is
interesting to place it in the context of other work in the
literature.  Let $T(Y)$ denote the set $T(\ma{Y};\ma{a},H)$ in the
special case $Y_i=Y$ and $a_i=H=1$, and let $T_*(Y)$
denote the set of $\ma{y}\in T(Y)$ for which the product $y_1y_2y_3$
is square-free.  Then it follows from Proposition \ref{guo} that
\beq\lab{pill}
\sum_{\ma{y}\in T_*(Y)} 2^{\omega(y_1y_2y_3)}
\leq \sum_{\ma{y}\in T(Y)} 2^{\omega(y_1y_2y_3)}
\ll Y^3.
\eeq
This should be compared with the work of Guo \cite{g} and Serre
\cite{serre}, who have both made a study of the corresponding sum weighted
instead by $1$.  By using the large sieve inequality Serre has shown that 
$$
\#T(Y)\ll \frac{Y^3}{(\log Y)^{3/2}},
$$
and Guo has proved an asymptotic formula for $\#T_*(Y)$ which
agrees with this upper bound.  In particular,
together with (\ref{pill}),  these estimates show that
a random conic in $\bfP^2$ does not contain a rational point.

We shall need several technical results during the proof of Proposition
\ref{guo}, and it will be convenient to list them here.  We begin by 
recording a basic estimate for the average orders of $\vt$ and 
$\vt^2$, as given by (\ref{phi}), whose simple proof we include here for completeness.

\begin{lem}\lab{vt-up}
Let $z\geq 1$.  Then we have
$$
\sum_{n\leq z} \vartheta(n) \leq \sum_{n\leq z} \vartheta(n)^2 \ll z.
$$
\end{lem}

\begin{proof}
The first inequality is trivial, since $\vt(n) \geq 1$ for all $n \in
\N$.  To see the second inequality we note that
\beq\lab{val}
\sum_{n\leq z} \vartheta(n)^2 =\sum_{n \leq z}\Big(\sum_{d\mid
n}\frac{|\mu(d)|}{d}\Big)^2
=\sum_{n \leq z}\sum_{[d_1,d_2]\mid n}\frac{|\mu(d_1)\mu(d_2)|}{d_1d_2},
\eeq
where $[d_1,d_2]=d_1d_2/\hcf(d_1,d_2)$ denotes the least
common multiple of $d_1$ and $d_2$.  But on writing $n=[d_1,d_2]e$ we
easily  deduce that
$$
\sum_{n\leq z} \vartheta(n)^2 \leq \sum_{d_1,d_2 =1}^\infty\sum_{e \leq
z/[d_1,d_2]}\frac{1}{d_1d_2} \leq
z\sum_{d_1,d_2=1}^\infty\frac{\hcf(d_1,d_2)}{d_1^2d_2^2}\ll z,
$$
since
$$
\sum_{d_1,d_2=1}^\infty\frac{\hcf(d_1,d_2)}{d_1^2d_2^2}\leq
\sum_{d_1,d_2=1}^\infty\sum_{k \mid d_1,d_2}\frac{k}{d_1^2d_2^2}\leq
\sum_{k,d_1',d_2'=1}^\infty\frac{1}{k^3{d_1'}^2{d_2'}^2}\ll 1.
$$
This completes the proof of Lemma \ref{vt-up}.
\end{proof}

We shall also need to apply the large sieve inequality in our work.
The following version of the large sieve is due to Montgomery  \cite{large}.

\begin{lem}\lab{sieve}
Let $z, N \geq 1$.  Suppose that 
$S\subseteq \Z\cap [1,N],$ and that for every prime $p \leq z$ there exists $\sigma(p) \in [0,p)$ 
such that the image of $S$ in $\Z/p\Z$ has $p-\sigma(p)$
elements.   Then we have
$$
\#S \ll \frac{N+z^2}{G(z)},
$$
with
$$
G(z)=\sum_{n \leq z}|\mu(n)|\prod_{p \mid n} \frac{\sigma(p)}{p-\sigma(p)}. 
$$
\end{lem}

During the course of this section we will need some standard character
sum estimates.  
The following estimate is due to Heath-Brown
\cite[Corollary 4]{h-b95}.

\begin{lem}\lab{hb95}
Let $M,N \in \N$, and let $a_1,\ldots,a_M$ and $b_1,\ldots,b_N$ be
arbitrary complex numbers satisfying $|a_m|,|b_n|\leq 1$.  Then 
$$
\sum_{\colt{m \leq M}{2\nmid m}}\sum_{n \leq N}a_mb_n
\Big(\frac{n}{m}\Big) \ll_\ve (MN)^\ve (M^{1/2}N+MN^{1/2}),
$$
for
any $\ve>0$.
\end{lem}

Next we recall the P\'olya--Vinogradov inequality, which 
can be found in the work of Davenport \cite[Chapter 23]{dav}, for example.

\begin{lem}\lab{pv}
Let $M,N \in \N$, and let $\chi$ be a non-principal character modulo
$q$.  Then
$$
\sum_{M \leq n \leq N} \chi(n) \ll q^{1/2}\log q.
$$
\end{lem}

We are now ready to commence our proof of Proposition \ref{guo} in
earnest.  It is important to note that in the course of this proof, as throughout our work,
we shall follow common practice and allow the small positive constant
$\ve$ to take different values at different parts of the argument.
For any $\ma{a}\in Z_*^3$, and any $Y_1,Y_2,Y_3 \geq 1$, recall the definition of the set
$T(\ma{Y};\ma{a},H)$ defined above.  It will be convenient to write
\beq\lab{defT}
\mcal{T}=\mcal{T}(\ma{Y};\ma{a},H)=\sum_{\ma{y}\in T(\ma{Y};\ma{a},H)} 2^{\omega(y_1y_2y_3)},
\eeq
with which notation our task is to establish the inequality 
$$
\mcal{T}\ll_\ve 
\vt(a_1a_2)H^\ve \Big(Y_1Y_2Y_3
+(Y_1Y_2)^{1/2+\ve}Y_3m_\ve(\ma{a,Y})
\Big),
$$
where $m_\ve(\ma{a,Y})$ is given by (\ref{m}).
Our approach is based upon a fusion of the ideas used by 
Guo and Serre in their work on this problem.  Recall the assumption (\ref{scat}).
Whenever $Y_2$ is much smaller in size than $Y_3$ we shall be in a
position to apply Lemma \ref{sieve} to estimate $\mcal{T}$.
In the alternative case, in which $Y_3$ is bounded above by a power of
$Y_2$, we shall employ Guo's approach.    We begin by studying the
latter case.

\subsubsection{Proof of Proposition \ref{guo}: $Y_3 \leq Y_2^{10}$}

As indicated by the title, our aim in this section is to establish
Proposition \ref{guo} under the assumption that
\beq\lab{assume-1}
Y_3 \leq Y_2^{10}.
\eeq
Given any $n \in \Z$, it will be convenient to extend the definition of
the Legendre symbol $(\frac{n}{p})$ to all primes $p$ by setting $(\frac{n}{2})=0$.
Our first step is the observation that
\beq\lab{legend}
1+\Big(\frac{n}{p}\Big)= \left\{
\begin{array}{ll}
2, & (\frac{n}{p})=1,\\
0, & (\frac{n}{p})=-1,\\
1, & \mbox{otherwise},\\
\end{array}
\right.
\eeq
for any integer $n$ and prime $p$.  Recalling (\ref{jac-1}) and
the fact that $\hcf(y_i,y_j)=1$ for any $\ma{y} \in
T(\ma{Y};\ma{a},H)$, we see that
\begin{align*}
\mcal{T}&=
\sum_{\ma{y}\in T(\ma{Y};\ma{a},H)}2^{\omega(y_1)+\omega(y_2)+\omega(y_3)}\\
&\ll_\ve H^\ve \sum_{\colt{|y_i| \leq Y_i}{\hcf(y_i,y_j)=1}}
\prod_{\{i,j,k\}=\{1,2,3\}}\prod_{p \mid y_k}\Big(1+\Big(\frac{-a_ia_jy_iy_j}{p}\Big)\Big)\\
&\ll_\ve H^\ve
\sum_{\colt{|y_i| \leq Y_i}{\hcf(y_i,y_j)=1}}
{\sum_{d_i \mid y_i}}^{\sharp} |\mu(d_1d_2d_3)|
\prod_{\{i,j,k\}=\{1,2,3\}}\Big(\frac{-a_ia_jy_iy_j}{d_k}\Big),
\end{align*}
where $\Sigma^\sharp$ denotes a summation over odd divisors $d_i \mid y_i$.
Here we have used (\ref{divH}) to deduce that
$2^{\omega(\hcf(y_k,a_ia_j))} \leq 2^{\omega(H)} \ll_\ve H^\ve$.
We first show that the contribution arising from case in
which $d_2=d_3=1$ is $O_\ve(H^\ve Y_1Y_2Y_3)$, which is satisfactory.
But the contribution from this case is clearly 
\begin{align*}
&\ll_\ve H^\ve Y_1Y_2Y_3 +H^\ve
\sum_{\colt{|y_i| \leq Y_i}{\hcf(y_i,y_j)=1}}
{\sum_{\colt{d_1 \mid y_1}{d_1 \neq 1}}}^\sharp |\mu(d_1)|
\Big(\frac{-a_2a_3y_2y_3}{d_1}\Big).
\end{align*}
Using the M\"obius function to pick out the condition
$\hcf(y_3,y_1y_1)=1$, we may therefore apply Lemma \ref{pv} to deduce that
\begin{align*}
&\ll_\ve H^\ve 
\sum_{|y_2|\leq Y_2} \sum_{\colt{|d_1e_1| \leq Y_1}{d_1 \neq 1}}
|\mu(d_1)|
\Big|\sum_{k \mid d_1e_1y_2}\sum_{|z_3| \leq Y_3/k} \mu(k)\Big(\frac{kz_3}{d_1}\Big)\Big|\\
&\ll_\ve H^\ve 
\sum_{|y_2|\leq Y_2} \sum_{\colt{|d_1e_1| \leq Y_1}{d_1 \neq 1}}
|\mu(d_1)|
\sum_{k\mid d_1e_1y_2}|\mu(k)|
\Big|\sum_{|z_3| \leq Y_3/k} \Big(\frac{z_3}{d_1}\Big)\Big|\\
&\ll_\ve (HY_1Y_2)^\ve 
\sum_{|y_2|\leq Y_2} \sum_{\colt{|d_1e_1| \leq Y_1}{d_1 \neq 1}}d_1^{1/2}\\
&\ll_\ve (HY_1Y_2)^\ve Y_1^{3/2}Y_2 \ll_\ve H^\ve Y_1Y_2Y_3,
\end{align*}
if $\ve>0$ is sufficiently small.
It therefore suffices to establish the estimate
\beq\lab{way1}
\sum_{\colt{|y_i| \leq Y_i}{\hcf(y_i,y_j)=1}}
{\sum_{\colt{d_i \mid y_i}{d_2d_3\neq 1}}}^\sharp |\mu(d_1d_2d_3)|
\prod_{\{i,j,k\}=\{1,2,3\}}\Big(\frac{-a_ia_jy_iy_j}{d_k}\Big) \ll Y_1Y_2Y_3,
\eeq
in order to complete the proof of Proposition \ref{guo}
under the assumption (\ref{assume-1}).  

Our main tool in the proof of (\ref{way1}) will be Lemma \ref{hb95}.
Let $D_1,D_2,D_3 \geq 1$.  We begin by estimating the contribution to the left hand side of
(\ref{way1}) from those values of $d_i$ such that
$$
D_i \leq d_i <2D_i.
$$
Let us write $\mcal{T}(\ma{D})$ for this contribution.  
Ultimately we shall sum over dyadic intervals for $D_i\ll Y_i$ to deduce
(\ref{way1}).  Now for any permutation $\{i,j,k\}$ of $\{1,2,3\}$ we
see that
\begin{align*}
\mcal{T}(\ma{D})&\ll \sum_{|d_ie_i|\leq Y_i}\sum_{d_j \ll
  D_j}\sum_{e_k \ll Y_k/D_k}\Big|  
{\sum_{D_k \leq d_k <2D_k}}^{\sharp} \sum_{e_j \ll Y_j/D_j} c(d_k,e_j) \Big(\frac{e_j}{d_k} \Big)
\Big|,
\end{align*}
for certain coefficients $c(d_k,e_j)\in \Z$ having modulus at most $1$.
Hence it easily follows from Lemma \ref{hb95} that for any $\ve>0$ we have
\begin{align}
\mcal{T}(\ma{D})
&\ll_\ve Y_3^\ve\sum_{|d_ie_i|\leq Y_i}\sum_{d_j \ll
  D_j}\sum_{e_k \ll
  Y_k/D_k}\Big(D_k^{1/2}\frac{Y_j}{D_j}+D_k\Big(\frac{Y_j}{D_j}\Big)^{1/2} \Big)\nonumber\\
&\ll_\ve Y_3^\ve\Big( \frac{Y_iY_jY_k}{D_k^{1/2}}+ Y_iY_j^{1/2}Y_kD_j^{1/2} \Big).\lab{km2}
 \end{align}
This will be satisfactory when exactly one of $D_2$ or $D_3$ is large.  
To handle the case in which both $D_2$ and $D_3$ are large, we proceed
by noting that 
\begin{align*}
\mcal{T}(\ma{D})&\ll \sum_{|d_1e_1|\leq Y_1}\sum_{\colt{e_2 \ll
  Y_2/D_2}{e_3 \ll Y_3/D_3}}\Big|  
{\sum_{d_2,d_3}}^{\sharp} c'(d_2,d_3) \Big(\frac{d_2}{d_3} \Big)\Big|,
\end{align*}
for certain coefficients $c'(d_2,d_3)\in \Z$ having modulus at most $1$.
But then an application of Lemma \ref{hb95} immediately yields
\beq\lab{km3}
\mcal{T}(\ma{D})\ll_\ve \frac{Y_1Y_2Y_3^{1+\ve}}{D_2D_3}\Big(D_2^{1/2}D_3+D_2D_3^{1/2}\Big),
\eeq
for any $\ve>0$.  
Finally we must consider the case in which both $D_2$ and $D_3$ are
small.  For this we recall that we have already handled the
contribution from those $d_2,d_3$ for which $d_2d_3=1$.  Suppose first that
$d_2\neq 1$. Then Lemma \ref{pv} yields the contribution
\begin{align*}
&\ll \sum_{d_3\ll D_3}\sum_{\colt{e_1 \ll Y_1/D_1}{e_2
    \ll Y_2/D_2}}\sum_{\colt{d_1 \ll D_1, d_2\ll D_2}{d_2 \neq 1}}|\mu(d_1d_2)|
\Big|  
\sum_{k \mid d_1d_2e_1e_2}\mu(k)\sum_{f_3 \ll Y_3/(D_3k)} \Big(\frac{kf_3}{d_1d_2} \Big)\Big|\\
&\ll_\ve Y_2^{\ve}
\sum_{d_3\ll D_3}\sum_{\colt{e_1 \ll Y_1/D_1}{e_2
    \ll Y_2/D_2}}\sum_{d_1 \ll D_1, d_2\ll D_2}
(D_1D_2)^{1/2}\\
&\ll_\ve Y_1Y_2^{1+\ve} (D_1D_2)^{1/2}D_3 \ll_\ve Y_1Y_2Y_3^{1/2+\ve} D_2^{1/2}D_3
\end{align*}
to $\mcal{T}(\ma{D})$, since $D_1 \ll Y_1 \leq Y_3$.  
Alternatively, if $d_2=1$ and $d_3\neq 1$, then we obtain the contribution
\begin{align*}
&\ll \sum_{\colt{e_1 \ll Y_1/D_1}{e_3
    \ll Y_3/D_3}}\sum_{\colt{d_1 \ll D_1, d_3\ll D_3}{d_3 \neq 1}}|\mu(d_1d_3)|
\Big|  \sum_{k \mid d_1d_3e_1e_3}\mu(k)
\sum_{f_2 \ll Y_2/k} \Big(\frac{kf_2}{d_1d_3} \Big)\Big|\\
&\ll_\ve Y_1Y_2^{1/2}Y_3^{1+\ve} D_3^{1/2}
\end{align*}
to $\mcal{T}(\ma{D})$.  Thus we may combine these two estimates to
deduce that
\beq\lab{km1}
\mcal{T}(\ma{D}) \ll_\ve Y_1Y_2^{1/2}Y_3^{1/2+\ve}
\Big(D_3^{1/2}Y_3^{1/2}+Y_2^{1/2}D_2^{1/2}D_3\Big). 
\eeq

We are now in a position to collect together our various estimates to
establish the bound
\beq\lab{way3}
\mcal{T}(\ma{D}) \ll_\ve Y_1Y_2^{7/8}Y_3^{1+\ve},
\eeq
for any $\ve>0$.  Before doing so we note that on summing over dyadic
intervals for $D_i \ll Y_i$ this is enough to establish that the left
hand side of (\ref{way1}) is  $O_\ve (Y_1Y_2^{7/8}Y_3^{1+\ve})$.  It is
at this point that we employ the assumption (\ref{assume-1}), which implies in
particular that $\log Y_3 \ll \log Y_2$.  This therefore establishes that
$$
\mcal{T} \ll_\ve H^\ve Y_1Y_2Y_3
$$ 
provided that the value of $\ve$ is taken to
be sufficiently small, and so completes the deduction of Proposition
\ref{guo} from (\ref{way3}) under the assumption that (\ref{assume-1})
holds. 
In order to establish (\ref{way3}) we shall need to split the argument
according to the size of $D_2, D_3$.  On supposing first that $D_2,D_3 \ll
Y_2^{1/4}$, it clearly follows from (\ref{km1}) that
$$
\mcal{T}(\ma{D}) \ll_\ve Y_1Y_2^{1/2}Y_3^{1/2+\ve}
\Big(Y_2^{1/8}Y_3^{1/2}+Y_2^{7/8}\Big) \ll_\ve
Y_1Y_2^{7/8}Y_3^{1+\ve},
$$
which is satisfactory for (\ref{way3}). Similarly, if 
$D_2,D_3 \gg Y_2^{1/4}$, then (\ref{km3}) implies that
(\ref{way3}) holds.  Next we suppose that  $D_2 \ll 
Y_2^{1/4}$ and $D_3 \gg Y_2^{1/4}$.  But then we may apply
(\ref{km2}) to the permutation $(i,j,k)=(1,2,3)$ to get
$$
\mcal{T}(\ma{D}) \ll_\ve 
Y_3^\ve\Big(\frac{Y_1Y_2Y_3}{Y_2^{1/8}}+ Y_1Y_2^{5/8}Y_3\Big)
\ll_\ve Y_1Y_2^{7/8}Y_3^{1+\ve},
$$
which is satisfactory for (\ref{way3}).
Finally, if  $D_2 \gg Y_2^{1/4}$ and $D_3 \ll Y_2^{1/4}$ then an
application of (\ref{km2}) to the permutation $(i,j,k)=(1,3,2)$ also
yields (\ref{way3}).

\subsubsection{Proof of Proposition \ref{guo}: $Y_3 > Y_2^{10}$}

We now turn to the proof of Proposition \ref{guo} under the assumption that
\beq\lab{assume-2}
Y_3 > Y_2^{10}.
\eeq
In view of the previous section, this will suffice to complete the
proof of Proposition \ref{guo}.
Note that $2^{\omega(n)} \leq d(n)$ for any $n \in \N$, where $d(n)$ is the ordinary divisor
function.  Then on recalling the conditions (\ref{divH})  and (\ref{jac-1}), 
we see that the quantity (\ref{defT}) satisfies
\begin{align*}
\mcal{T}
&\leq \sum_{\colt{|y_1| \leq Y_1, |y_2|\leq
    Y_2}{\hcf(y_1,y_2)=1}}2^{\omega(y_1y_2)}
\sum_{y_3: ~\ma{y}\in T(\ma{Y};\ma{a},H)}\sum_{d \mid y_3}1\\
&\ll \sum_{\colt{|y_1| \leq Y_1, |y_2|\leq
    Y_2}{\hcf(y_1,y_2)=1}}2^{\omega(y_1y_2)}
{\sum_{d\leq Y_3^{1/2}}}^{\dagger}
\#\{e \in \Z: (y_1,y_2,de) \in T(\ma{Y};\ma{a},H)\},
\end{align*}
where the summation $\Sigma^\dagger$ is only over integers $d\leq
Y_3^{1/2}$ such that 
\beq\lab{dj}
\hcf(d,a_1a_2)\mid H^2, \quad \hcf(d,y_1y_2)=1, 
\eeq
and
$$
\Big(\frac{-a_1a_2y_1y_2}{p}\Big)=1
$$
for all odd primes $p$ dividing $d/\hcf(d,a_1a_2)$.

On defining the set
$$
S=S(Y_3/d;y_1,y_2,\ma{a},H)=\{e \in \Z: (y_1,y_2,de) \in T(\ma{Y};\ma{a},H)\},
$$
so that in particular 
\beq\lab{tache}
\mcal{T}
\ll \sum_{\colt{|y_1| \leq Y_1, |y_2|\leq
    Y_2}{\hcf(y_1,y_2)=1}}2^{\omega(y_1y_2)}
{\sum_{d\leq Y_3^{1/2}}}^{\dagger} \#S,
\eeq
we see that our task is now to estimate the size of $S$.  For this 
we shall use the large sieve inequality, as presented in Lemma \ref{sieve}.
For any odd prime $p$ we must calculate the size of the image $S_p$
of $S$ in $\Z/p\Z$.   
If $p \nmid a_iy_1y_2d$ then the congruence
$$
a_1y_1x_1^2+a_2y_2x_2^2+a_3dex_3^2\eqm{0}{p}
$$
is always soluble whenever $p \nmid e$, by the
Chevalley--Warning theorem. Alternatively, if $p\mid e$ then this
congruence is soluble if and only if
\beq\lab{payne}
\Big(\frac{-a_1a_2y_1y_2}{p}\Big)=1.
\eeq
Hence we conclude that 
$$
\#S_p= 
\left\{
\begin{array}{ll}
p-1, & (\frac{-a_1a_2y_1y_2}{p})=-1,\\
p, & \mbox{otherwise},
\end{array}
\right.
$$
whenever $p \nmid a_iy_1y_2d$.
Turning to the size of $S_p$ in the case $p \mid a_iy_1y_2$, we
suppose that $p \mid a_1y_1$ and $p \nmid a_2a_3y_2d$.  But then either $e=0$, or else $e \neq
0$ and $e$ belongs to $S_p$ if and only if 
$$
\Big(\frac{-a_2a_3y_2de}{p}\Big)=1.
$$
Hence there are $\frac{1}{2}(p+1)$
possible values of $e$ over all.  Similarly one finds that
$\#S_p=\frac{1}{2}(p+1)$ if $p \mid a_2y_2$ and $p \nmid a_1a_3y_1d$.  
If $p \mid a_3d$ and $p \nmid a_1a_2y_1y_2$ however, then $\#S_p=p$
since we may clearly suppose that (\ref{payne}) holds automatically for such primes.
Finally we note that $\#S_p=p$ in all other cases.
Taking $\sigma(p)$ to be $p-\#S_p$, we have therefore shown that
$$
\frac{\sigma(p)}{p-\sigma(p)} = 
\left\{
\begin{array}{ll}
1/(p-1), & p\nmid a_iy_1y_2d, (\frac{-a_1a_2y_1y_2}{p})=-1,\\
(p-1)/(p+1), & p \mid a_1y_1, p\nmid a_2a_3y_2d,\\
(p-1)/(p+1), & p \mid a_2y_2, p\nmid a_1a_3y_1d,\\
0, & \mbox{otherwise}.
\end{array}
\right.
$$
Now write $g(n)=|\mu(n)|\prod_{p \mid n} \frac{\sigma(p)}{p-\sigma(p)}$ for
any $n \in \N$, so that in particular
$$
g(mn)=\left\{
\begin{array}{ll}
g(m)g(n), & \mbox{if $\hcf(m,n)=1$,}\\
0, & \mbox{otherwise,}
\end{array}
\right.
$$
and $G(z)=\sum_{n \leq z}g(n)$.  But then it is easy to deduce that 
$$
G(z)= \sum_{k \mid a_iy_1y_2d}g(k)\sum_{\colt{n\leq
    z/k}{\hcf(n,a_iy_1y_2d)=1}}g(n) \geq \sum_{k \mid y_1y_2}g(k)\sum_{\colt{n\leq z/k}{\hcf(n,a_iy_1y_2d)=1}}g(n),
$$
for any $z \geq 1$.  We may now use the coprimality condition
$\hcf(y_1,y_2)=1$, together with (\ref{divH}) and (\ref{dj}) to conclude that 
$$
\sum_{k \mid y_1y_2}g(k)=\prod_{p \mid y_1y_2}\Big(1+g(p)\Big)
\gg_\ve H^{-\ve}2^{\omega(y_1y_2)}\prod_{p \mid y_1y_2}\Big(1+\frac{1}{p}\Big)^{-1},
$$
for any $\ve>0$.  Moreover it is not hard to deduce that
\begin{align*}
\sum_{\colt{n\leq z/k}{\hcf(n,a_iy_1y_2d)=1}}g(n)
&\gg \sum_{\colt{n\leq z/k}{p\mid n \Rightarrow
    (\frac{-a_1a_2y_1y_2}{p})=-1}}\frac{|\mu(n)|}{n}\\
&\gg \log(z/k) \Big(\sum_{\colt{n\leq z/k}{p\mid n \Rightarrow
    (\frac{-a_1a_2y_1y_2}{p})\neq -1}}\frac{|\mu(n)|}{n}\Big)^{-1}.
\end{align*}
We may now complete our estimate for $\#S$ by taking $z=Y_3^{1/4}$ in
Lemma \ref{sieve}.  On noting that $Y_3/d \geq Y_3^{1/2}$, since
$d \leq Y_3^{1/2}$, it therefore follows that
$$
\#S \ll \frac{Y_3}{dG(Y_3^{1/4})}.
$$
It is at this point that we apply the hypothesis (\ref{assume-2}), from
which it follows that
$$
Y_1Y_2 \leq Y_2^2 <Y_3^{1/5}.
$$
In particular we see that $Y_3^{1/4}/k> Y_3^{1/20}$
for any divisor $k \mid y_1y_2$.  Recall the definition (\ref{phi}) of the arithmetic function
$\vt$.   We therefore deduce that
$$
G(Y_3^{1/4}) \gg_\ve \frac{2^{\omega(y_1y_2)}}{H^\ve\vt(y_1y_2)}
\log Y_3
\Big(\sum_{\colt{n\leq Y_3^{1/2}}{p\mid n \Rightarrow
    (\frac{-a_1a_2y_1y_2}{p})\neq -1}}\frac{|\mu(n)|}{n}\Big)^{-1},
$$
for any $\ve>0$, whence
$$
\#S \ll_\ve \frac{H^\ve \vt(a_1a_2)\vt(y_1y_2)^2}{2^{\omega(y_1y_2)}} \frac{Y_3}{d\log Y_3}
\sum_{\colt{n\leq Y_3^{1/2}}{p\mid n \Rightarrow
    (\frac{-a_1a_2y_1y_2}{p})=1}}\frac{1}{n}.
$$
On inserting this into (\ref{tache}), we therefore deduce the statement
of the following result.

\begin{lem}\lab{k=1}
Let $\ve>0$.  Then we have
$$
\mcal{T}
\ll_\ve \vt(a_1a_2)H^\ve\frac{Y_3}{\log Y_3}
\sum_{\colt{|y_1| \leq Y_1, |y_2|\leq
    Y_2}{\hcf(y_1,y_2)=1}}
\vt(y_1y_2)^2 
\Big(\sum_{n \in \mcal{N}}
\frac{1}{n}\Big)^2,
$$
with
$$
\mcal{N}=\mcal{N}(Y_3;\ma{a},y_1,y_2)=\Big\{n \in \N: n\leq Y_3^{1/2},
\;\mbox{$\Big(\frac{-a_1a_2y_1y_2}{p}\Big)=1$ for all odd $p \mid n$}\Big\}.
$$ 
\end{lem}

On writing $\mcal{M}=\mcal{N}(Y_3^2;\ma{a},y_1,y_2)$, 
we observe that
$$
\Big(\sum_{n \in \mcal{N}}
\frac{1}{n}\Big)^2
= \sum_{n_1, n_2 \in \mcal{N}}
\frac{1}{n_1n_2}\leq \sum_{m \in \mcal{M}}\frac{d(m)}{m},
$$
in the statement of Lemma \ref{k=1}.  But now we  may clearly apply
(\ref{legend}) in Lemma~\ref{k=1} to deduce that
\beq\lab{chav}
\mcal{T}
\ll_\ve \frac{\vt(a_1a_2)H^\ve Y_3}{\log Y_3}
\sum_{\colt{|y_1| \leq Y_1, |y_2|\leq
    Y_2}{\hcf(y_1,y_2)=1}}
\vt(y_1y_2)^2
\sum_{m \leq Y_3}\frac{1}{m}
\sum_{d \mid m}\Big(\frac{-a_1a_2y_1y_2}{d}\Big).
\eeq

We proceed by considering the
contribution to (\ref{chav}) from the terms for which $d=k^2$ is a square.  Writing
$m=jk^2$, we obtain the contribution
\begin{align*}
&\ll_\ve
\frac{\vt(a_1a_2)H^\ve Y_3}{\log Y_3}
\sum_{y_1,y_2} \vt(y_1)^2\vt(y_2)^2
\sum_{k \leq Y_3^{1/2}}\sum_{j \leq Y_3/k^2}\frac{1}{jk^2}\\
&\ll_\ve 
\frac{\vt(a_1a_2)H^\ve Y_3}{\log Y_3}
\sum_{j \leq Y_3}\frac{1}{j}
\sum_{y_1,y_2} \vt(y_1)^2\vt(y_2)^2\\
&\ll_\ve  \vt(a_1a_2)H^\ve Y_1Y_2Y_3,
\end{align*}
by Lemma \ref{vt-up}.  This is plainly satisfactory for Proposition
\ref{guo}. 
In order to handle the contribution to (\ref{chav}) from the remaining divisors, we
define the characteristic function 
$$
\delta(n)=
\left\{
\begin{array}{ll}
0, & \mbox{$n=k^2$ for some $k \in \N$,}\\
1, & \mbox{otherwise,}
\end{array}
\right.
$$
for any $n \in \N$.
Note that in particular we have $\delta(1)=0$.
Writing $m=de$, it follows that 
\beq\lab{not1-x}
\sum_{\colt{|y_1| \leq Y_1, |y_2|\leq
    Y_2}{\hcf(y_1,y_2)=1}}
\vt(y_1y_2)^2
\sum_{m \leq Y_3}\frac{1}{m}
\sum_{\colt{d \mid m}{\delta(d)=1}}\Big(\frac{-a_1a_2y_1y_2}{d}\Big)
=
\sum_{e \leq Y_3}\frac{S_{e}}{e},
\eeq
with
$$
S_{e}= 
\sum_{d \leq Y_3/e} \frac{\delta(d)}{d} \sum_{\colt{|y_1| \leq Y_1, |y_2|\leq
    Y_2}{\hcf(y_1,y_2)=1}}
\vt(y_1)^2\vt(y_2)^2\Big(\frac{-a_1a_2y_1y_2}{d}\Big).
$$
Our next task is to establish the following inequality.

\begin{lem}\lab{sex}
Let $\ve>0$. Then we have
$$
S_{e}\ll_{\ve}
Y_1Y_2+(Y_1Y_2)^{1/2+\ve}(\min\{|a_1a_2|,Y_3\}^{\ve}+\log Y_3).
$$
\end{lem}

\begin{proof}
Let us consider the contribution $S_{e}(D_1,D_2)$ to $S_{e}$ from $d$ contained in the
interval $D_1\ll d \ll D_2$, for various choices of  
\beq\lab{burg-x}
1\leq D_1 \leq D_2 \leq Y_3/e.
\eeq
Suppose that $N \in \Z$ is not a square, and that $D \geq 1$.  Then
an application of Lemma \ref{pv} yields
$$
\sum_{d \leq D}\delta(d) \Big(\frac{N}{d}\Big)= 
\sum_{d \leq D}\Big(\frac{N}{d}\Big)- \sum_{\colt{d \leq D}{d=k^2}}1
\ll N^{1/2}\log N +D^{1/2}.
$$
Thus for any $D_1, D_2$ in the range (\ref{burg-x}), we may combine partial
summation with Lemma \ref{vt-up} to deduce that the contribution to 
$S_{e}(D_1,D_2)$ from those $y_1,y_2$ for which $-a_1a_2y_1y_2$
is not a square, is 
\begin{align*}
&\ll \sum_{y_1,y_2}\delta(-a_1a_2y_1y_2)\vt(y_1)^2\vt(y_2)^2\Big|\sum_{D_1 \ll d \ll D_2}
  \frac{\delta(d)}{d} \Big(\frac{-a_1a_2y_1y_2}{d}\Big) \Big|\\
&\ll_\ve \sum_{y_1,y_2} \vt(y_1)^2\vt(y_2)^2\Big(D_1^{-1/2} + D_1^{-1}|a_1a_2y_1y_2|^{1/2+\ve}\Big)\\
&\ll_\ve Y_1Y_2 +
  D_1^{-1}|a_1a_2|^{1/2+\ve}(Y_1Y_2)^{3/2+\ve}\\
&\ll Y_1Y_2 + D_1^{-1}|a_1a_2|^{3/4}(Y_1Y_2)^{7/4},
\end{align*}
on taking $\ve=1/4$.  Now there are clearly $O((Y_1Y_2)^{1/2})$ values of $y_1,y_2$ for which
$-a_1a_2y_1y_2$ is a square.  
In view of the trivial inequality $\vt(n)\ll_{\ve}n^{\ve}$ for any $\ve>0$, it follows
that the total contribution to $S_{e}(D_1,D_2)$ from such $y_1,y_2$ is
$$
\ll_{\ve}(Y_1Y_2)^{1/2+\ve} \log D_2
$$ 
for any $\ve>0$. We have therefore established the bound
\beq\lab{cd3-x}
S_e(D_1,D_2) \ll_{\ve} Y_1Y_2 + D_1^{-1}|a_1a_2|^{3/4}(Y_1Y_2)^{7/4}
+(Y_1Y_2)^{1/2+\ve} \log D_2,
\eeq
for any  $\ve>0$ and any $D_1, D_2$ in the range (\ref{burg-x}).

We shall need an alternative estimate for $S_e(D_1,D_2)$ to handle the contribution
from small values of $d$. 
Let $[d_1,d_2]$ denote the least common multiple of $d_1$ and $d_2$, as usual.
Then it follows
from a combination of Lemma \ref{pv} and (\ref{val}), that for fixed values of $d, y_1$ such that
$\delta(d)=1$, we have
\begin{align*}
\Big|\sum_{\colt{y_2\leq Y_2}{\hcf(y_1,y_2)=1}} \vt(y_2)^2\Big(\frac{y_2}{d}\Big)\Big|
&\leq
\sum_{k\mid y_1}|\mu(k_1)|\vt(k)^2\Big|\sum_{y_2\leq Y_2/k} \vt(y_2)^2\Big(\frac{y_2}{d}\Big)\Big|\\
&\leq 
\sum_{k\mid y_1}\vt(k)^2
\Big|\sum_{[d_1,d_2]|e| \leq
  Y_2/k} \frac{|\mu(d_1)\mu(d_2)|}{d_1d_2} \Big(\frac{[d_1,d_2]e}{d}\Big) \Big|\\
&\leq \sum_{k\mid y_1}\vt(k)^2
\sum_{d_1,d_2\leq Y_2}
\frac{1}{d_1d_2} \Big| \sum_{|e| \leq Y_2/(k[d_1,d_2])}  \Big(\frac{e}{d}\Big)\Big|\\
&\ll_\ve d^{1/2+\ve}(Y_1Y_2)^{\ve}.
\end{align*}
Hence a further application of Lemma \ref{vt-up} yields
\begin{align}
S_{e}(1,(Y_1Y_2)^{1/2})
&\ll \sum_{|y_1|\leq Y_1}\vt(y_1)^2\sum_{d
  \leq (Y_1Y_2)^{1/2}} \frac{\delta(d)}{d} \Big|
\sum_{\colt{y_2\leq Y_2}{\hcf(y_1,y_2)=1}} \vt(y_2)^2\Big(\frac{y_2}{d}\Big)\Big|
\nonumber\\
&\ll_\ve \sum_{|y_1|\leq Y_1}\vt(y_1)^2\sum_{d
  \leq (Y_1Y_2)^{1/2}} d^{-1/2+\ve} (Y_1Y_2)^{\ve}\nonumber\\
&\ll_\ve \sum_{|y_1|\leq Y_1}\vt(y_1)^2(Y_1Y_2)^{1/4+\ve}\nonumber\\
&\ll  Y_1Y_2,\lab{cd1-x}
\end{align}
provided that $\ve>0$ is taken to be sufficiently small.

For larger $d$, we employ Lemma \ref{hb95}.
Suppose that $(Y_1Y_2)^{1/2}\leq  D \leq Y_3/e$, 
and write $a(d)=\delta(d)(\frac{-a_1a_2}{d})$.
Then proceeding as above we see that 
$$
S_{e}(D,D)=  \sum_{D \ll d \ll D} \frac{a_d}{d}
\sum_{\colt{[d_1,e_1]|f_1| \leq
  Y_1}{[d_2,e_2]|f_2| \leq Y_2}} 
\frac{b(\ma{d,e,f})}{d_1d_2e_1e_2}
  \Big(\frac{[d_1,e_1][d_2,e_2]f_1f_2}{d}\Big),
$$
where $b(\ma{d,e,f})=|\mu(d_1)\mu(d_2)\mu(e_1)\mu(e_2)|$
if $\hcf([d_1,e_1]f_1, [d_2,e_2]f_2)=1$ and $b(\ma{d,e,f})=0$ otherwise.
In particular we observe that $|a(d)b(\ma{d,e,f})| \leq 1$, and so 
there exists a certain constant $c(d,f_1,f_2)$ of
modulus at most $1$ such that
$$
S_{e}(D,D) \ll  
\sum_{\colt{d_1,e_1 \leq Y_1}{d_2,e_2 \leq Y_2}} 
\frac{1}{d_1d_2e_1e_2}
\Big|
\sum_{D \ll d \ll D} \frac{1}{d}
\sum_{\colt{|f_1| \leq
  Y_1/[d_1,e_1]}{|f_2| \leq Y_2/[d_2,e_2]}} 
c(d,f_1,f_2)
\Big(\frac{f_1f_2}{d}\Big)\Big|.
$$
On combining the fact that $Y_1Y_2\leq  D^2$ with Lemma \ref{hb95} we
therefore deduce that
\begin{align*}
S_{e}(D,D)
&\ll_{\ve} \frac{(Y_1Y_2)^{\ve}}{D^{1-\ve}}
\sum_{\colt{d_1,e_1 \leq Y_1}{d_2,e_2 \leq Y_2}}
\frac{1}{d_1d_2e_1e_2}
\Big( \frac{D^{1/2}Y_1Y_2}{[d_1,e_1][d_2,e_2]}+ \frac{D
  (Y_1Y_2)^{1/2}}{([d_1,e_1][d_2,e_2])^{1/2}}  \Big)\\
&\ll_{\ve}    D^{\ve-1/2}Y_1Y_2+ D^{\ve} (Y_1Y_2)^{1/2},
\end{align*}
for any $\ve>0$.  Summing over dyadic
intervals for $(Y_1Y_2)^{1/2}\leq  D \leq D_2$, for any choice of $D_2 \leq Y_3/e$,
we therefore obtain the estimate
$$
S_{e}((Y_1Y_2)^{1/2},D_2) \ll_{\ve}  Y_1Y_2+ D_2^{\ve}(Y_1Y_2)^{1/2}.
$$
Once combined with (\ref{cd1-x}) this yields the bound
\beq\lab{cd2-x}
S_{e}(1,D_2) \ll_{\ve}  Y_1Y_2+ D_2^{\ve}(Y_1Y_2)^{1/2},
\eeq
for any $(Y_1Y_2)^{1/2} \leq D_2 \leq Y_3/e$.

Taking $D_2=Y_3/e$ in (\ref{cd2-x}) we obtain the estimate
\beq\lab{cole1-x}
S_e \ll_{\ve} Y_1Y_2+ (Y_1Y_2)^{1/2}Y_3^\ve,
\eeq
for any $\ve>0$.  In order to obtain a second estimate, and so
complete the proof of Lemma \ref{sex}, we apply (\ref{cd2-x}) with
$D_2=|a_1a_2|^{3/4}(Y_1Y_2)^{3/4}$ and (\ref{cd3-x}) with
$D_1=|a_1a_2|^{3/4}(Y_1Y_2)^{3/4}$ and $D_2=Y_3/e$.  This produces the estimate
$$
S_e \ll_{\ve} 
Y_1Y_2+ (Y_1Y_2)^{1/2+\ve}\Big(|a_1a_2|^{\ve}+\log Y_3\Big),
$$ 
for any $\ve>0$.  
We complete the proof of Lemma \ref{sex} by taking the minimum of the
bounds provided by this inequality and (\ref{cole1-x}).
\end{proof}

It remains to substitute Lemma \ref{sex} into
(\ref{not1-x}), and then  insert the resulting estimate into
(\ref{chav}).  Thus it follows that
\begin{align*}
\mcal{T}
&\ll_\ve
\frac{\vt(a_1a_2)H^\ve Y_3}{\log Y_3}
\sum_{e \leq Y_3}\frac{S_{e}}{e}\\
&\ll_\ve
\vt(a_1a_2)H^\ve Y_3
\Big(Y_1Y_2+(Y_1Y_2)^{1/2+\ve}(\min\{|a_1a_2|,Y_3\}^{\ve}+\log Y_3)\Big),
\end{align*}
which thereby establishes Proposition \ref{guo} under the assumption
that (\ref{assume-2}) holds.

\section{The equation $a_1b_1c_1^2+a_2b_2c_2^2+a_3b_3c_3^2=0$}\lab{core-x}

The purpose of this section is to bring together the results in \S
\ref{diag} and \S \ref{sieve-in}, in order to make a study of the density of 
integer solutions to the equation 
\beq\lab{key-x}
a_1b_1c_1^2+a_2b_2c_2^2+a_3b_3c_3^2=0.
\eeq
Thus for any $A_i,B_i,C_i\geq 1$ we let 
$\mcal{M}(A_i,B_i,C_i)$ denote the number of $\ma{a,b,c} \in Z_*^3$
such that (\ref{key-x}) holds and
$$\lab{range-x}
|a_i| \leq A_i, \quad |b_i| \leq B_i, \quad |c_i| \leq C_i,
$$
with
\beq\lab{cat-x}
\hcf(a_i,c_j)=
\hcf(c_i,c_j)=1
\eeq
and 
\beq\lab{cut-x}
|\mu(a_1a_2a_3)|=1, \quad \hcf(a_i,b_j,b_k)=1.
\eeq
It will be convenient to set
$$
A=A_1A_2A_3, \quad B=B_1B_2B_3, \quad C=C_1C_2C_3.
$$
With this notation in mind, we proceed by establishing the following  result.

\begin{pro}\lab{M1-x}
For any $\ve>0$, we have
$$
\mcal{M}(A_i,B_i,C_i) \ll_\ve A^{2/3}B^{2/3}C^{1/3} +
\sigma\tau AB^{1/2}C^{1/2},
$$
where
$$
\sigma=1+\frac{\min\{A, B\}^\ve}{\min\{B_iB_j\}^{1/16}},
\quad 
\tau=1+\frac{\log B}{\min\{B_iB_j\}^{1/16}}.
$$
\end{pro}

\begin{proof}
We shall prove Proposition \ref{M1-x} by making suitable applications of Lemma~\ref{line}
and Lemma \ref{box}.  Our starting point is to fix choices of
$\ma{a,c} \in Z_*^3$, and count the 
corresponding number of $\ma{b} \in Z_*^3$  satisfying 
(\ref{key-x}) and $|b_i| \leq B_i$.
Applying Lemma~\ref{line} with
$$
\ma{h}=(a_1c_1^2,a_2c_2^2,a_3c_3^2),
$$ 
we easily obtain the upper bound
 \begin{align*}
\mcal{M}(A_i,B_i,C_i) &\ll \sum_{\ma{a,c}} \Big(1+
   \frac{B}{\max |a_i|c_i^2B_i}\Big) \\
 &\leq \sum_{\ma{a,c}} \Big(1+
B^{2/3}|a_1a_2a_3c_1^2c_2^2c_3^2|^{-1/3}\Big)\\
 &\ll AC +A^{2/3}B^{2/3}C^{1/3}.
 \end{align*}
We shall use this bound whenever $C \leq B$, under which hypothesis
the estimate in Proposition \ref{M1-x} clearly holds.

It remains to handle
the case in which  $C> B$.  For this we 
fix choices of $\ma{a,b} \in Z_*^3$ for which (\ref{cut-x}) holds, and count the corresponding
number of $\ma{c} \in Z_*^3$ satisfying (\ref{key-x}), (\ref{cat-x}) and $|c_i| \leq C_i$.
Thus we are in a position to apply Lemma~\ref{box} with 
$$
\ma{g}=(a_1,a_2,a_3), \quad \ma{h}=(b_1, b_2, b_3).
$$
In particular it follows from (\ref{key-x})---(\ref{cut-x})
that
\begin{align*}
D_{\ma{a,b}} &=
\hcf(b_1b_2,b_1b_3,b_2b_3)\hcf(a_1,b_2b_3)\hcf(a_2,b_1b_3)\hcf(a_3,b_1b_2)
\\
&=\hcf(b_1b_2,b_1b_3,b_2b_3)\\
&\leq \hcf(b_1,b_2)\hcf(b_1,b_3)\hcf(b_2,b_3)=E_\ma{b},
\end{align*}
say.
Moreover, since $|b_1b_2b_3| \leq B$ and $C>B$, we also have
$$
\frac{C}{|b_1b_2b_3|} \geq \frac{C}{B} > 1.
$$
Thus under the assumption $C>B$ we may conclude from Lemma
\ref{box} that
$$
\mcal{M}(A_i,B_i,C_i) \ll C^{1/2}{\sum_{\ma{a,b}}}^{*} \frac{E_\ma{b}^{3/4}}{|b_1b_2b_3|^{1/2}}
2^{\omega(b_1b_2b_3)},
$$
where $\sum_{\ma{a,b}}^{*}$ indicates a summation over 
$\ma{a, b}\in Z_*^3$ for which  $|a_i| \leq A_i$, $|b_i| \leq B_i$, (\ref{cut-x}) 
holds, and the equation (\ref{key-x}) has a solution $\ma{c} \in
Z_*^3$ with (\ref{cat-x}) holding.

In order to handle the term $E_\ma{b}^{3/4}$ in our estimate for $\mcal{M}(A_i,B_i,C_i)$, we write
$$
b_i=h_{ij}h_{ik}b_i', 
$$
for fixed $h_{12}, h_{13}, h_{23} \in \N$ such that $h_{ij}=h_{ji}$.
Then 
\beq\lab{Wi-x}
|b_i'| \leq \frac{B_i}{h_{ij}h_{ik}}=B_i',
\eeq
say.
Since $\ma{b}$ is primitive, it follows that 
$\hcf(h_{ij},h_{ik})=1$. Moreover, for fixed values of $h_{12},h_{13},h_{23}$, it
suffices to sum over $\ma{a, b'}\in Z_*^3$ for which
\beq\lab{black-x}
\hcf(h_{ij},a_kb_k')= \hcf(b_i',b_j')=1,
\eeq
by (\ref{cut-x}) and the fact that $\ma{b}$ is primitive.  With this
change of variables, the equation (\ref{key-x}) clearly becomes
\beq\lab{mouse2-x}
\be_1 b_1'c_1^2+ \be_2b_2'c_2^2+ \be_3b_3'c_3^2=0,
\eeq
where we have written 
$$
\be_i=a_ih_{ij}h_{ik}
$$ 
for fixed values of $a_i,h_{ij}$.
We shall need to record the equality
\beq\lab{secant-x}
\hcf(\be_ib_i',\be_jb_j')= h_{ij},
\eeq
which easily follows from combining the coprimality conditions
(\ref{cat-x}),(\ref{cut-x}), (\ref{black-x}) and
$\hcf(h_{ij},h_{ik})=1$, with the equation (\ref{mouse2-x}).

Write $h=h_{12}h_{13}h_{23}$ and suppose without loss of generality
that 
$$B_1' \leq B_2' \leq
B_3'.
$$  
Then upon collecting our work
together, we see that
$$
\mcal{M}(A_i,B_i,C_i) \ll C^{1/2}\sum_{h_{ij}}\sum_{\ma{a}}
\sum_{\ma{b'}\in \mcal{B}} \frac{2^{\omega(h^2b_1'b_2'b_3')}}{h^{1/4}|b_1'b_2'b_3'|^{1/2}},
$$
where $\mcal{B}=\mcal{B}(B_i;a_i,h_{ij})$ denotes the set of $\ma{b'} \in Z_*^3$
with pairwise coprime components, 
for which (\ref{Wi-x}) and (\ref{secant-x}) hold,  and (\ref{mouse2-x})
has a solution $\ma{c} \in Z_*^3$ with 
$\hcf(c_i,c_j)=1$.   
Using the trivial upper bound $2^{\omega(n)}\ll_\ve n^{\ve}$, we may
therefore combine Proposition \ref{guo} with Lemma \ref{vt-up} and partial summation to deduce that
\begin{align*}
\mcal{M}(A_i,B_i,C_i) &\ll_\ve C^{1/2}\sum_{h_{ij}}h^{\ve-1/4}\sum_{\ma{a}}
\sum_{\ma{b'}\in \mcal{B}}
\frac{2^{\omega(b_1'b_2'b_3')}}{|b_1'b_2'b_3'|^{1/2}}\\
&\ll_\ve AC^{1/2}\sum_{h_{ij}}\frac{(B_1'B_2'B_3')^{1/2}}{h^{1/4-\ve}}\Big(1
+\frac{\min\{A,B\}^{\ve}+\log B}{(B_1'B_2')^{1/2-\ve}}\Big),
\end{align*}
On recalling the definition (\ref{Wi-x}) of $B_i'$, we see that
$$
\sum_{h_{ij}}
\frac{(B_1'B_2'B_3')^{1/2}}{h^{1/4-\ve}} \ll B^{1/2},
$$
provided that $\ve>0$ is taken to be sufficiently small.
Similarly, using the inequalities 
$$
h_{12} \leq (B_1B_2)^{1/2}, \quad h_{13} \leq B_1/h_{12}, \quad h_{23}
\leq B_2/h_{12},
$$ 
we easily check that
\begin{align*}
\sum_{h_{ij}}h^{\ve-1/4}(B_1'B_2')^\ve{B_3'}^{1/2} &\ll
(B_1B_2)^{3/8+\ve}B_3^{1/2}\ll
(B_1B_2)^{7/16}B_3^{1/2},
\end{align*}
provided that $\ve>0$ is sufficiently small.
It therefore follows that the estimate in Proposition \ref{M1-x} holds in
the case $C>B$, and so holds unconditionally. 
\end{proof}

It turns out that we shall need an alternative estimate for
$\mcal{M}(A_i,B_i,C_i)$ to handle the case in which $B_1, B_2, B_3$ have particularly awkward
sizes.  The following result is rather easy to establish.

\begin{pro}\lab{M2-x}
We have
$$
\mcal{M}(A_i,B_i,C_i) \ll AB_iB_j(C_k +  C_iC_jA_k^{-1}) (\log AC)^2,
$$
for any permutation $\{i,j,k\}$ of the set $\{1,2,3\}$.  
\end{pro}

\begin{proof}
Our proof of Proposition \ref{M2-x} is based upon Heath-Brown's
treatment \cite[Lemma 4]{cayley} of the equation 
$
n_1^2n_2n_3+n_4^2n_5n_6=n_7n_8.
$
For fixed integers $a,b, q$ we let $\rho(q;a,b)$ denote the number of solutions to the
congruence 
$
at^2 + b \eqm{0}{q}.
$
For any value of $q$, we then have 
\beq\lab{broken}
\rho(q;a,b)\leq \sum_{d \mid q} |\mu(d)| \Big(\frac{-ab}{d}\Big).
\eeq
We shall establish Proposition \ref{M2-x} in the case $(i,j,k)=(1,2,3)$, say.
The other cases will follow by symmetry.
Now it follows from (\ref{key-x}) that for given $a_i, b_1,b_2,c_3$, and each corresponding
solution $t$ of the congruence
$$
a_1b_1t^2+a_2b_2\eqm{0}{a_3c_3^2},
$$
we must have
$c_1 \eqm{t c_2}{a_3c_3^2}.$    This gives rise to an equation of the
form $\ma{h.w}=0$, with $\ma{h}=(1,-t,a_3c_3^2)$ and
$\ma{w}=(c_1,c_2,k)$.  Upon recalling that
$\hcf(c_1,c_2)=1$ from (\ref{cat-x}), an application of Lemma
\ref{line} therefore yields the bound 
$$
\ll \rho(a_3c_3^2;a_1b_2,a_2b_2)\Big(1+ \frac{C_1C_2}{|a_3c_3^2|}\Big),
$$ 
for the number of possible $b_3, c_1,c_2$ given fixed choices for
$a_i,b_1,b_2$ and $c_3$.  On employing (\ref{broken}) we therefore have
\begin{align*}
\mcal{M}(A_i,B_i,C_i) 
&\ll 
\sum_{a_i,b_1,b_2,c_3} \rho(a_3c_3^2;a_1b_2,a_2b_2)\Big(1+ \frac{C_1C_2}{|a_3c_3^2|}\Big)\\ 
&\ll 
\sum_{a_i,b_1,b_2,c_3}\sum_{d
\mid a_3c_3} |\mu(d)| 
\Big(\frac{-a_1a_2b_1b_2}{d}\Big)\Big(1+ \frac{C_1C_2}{|a_3c_3^2|}\Big)\\
&\ll 
\sum_{a_i,b_1,b_2,c_3} d(a_3)d(c_3)\Big(1+ \frac{C_1C_2}{|a_3c_3^2|}\Big),
\end{align*}
since the sum over square-free divisors of $a_3c_3^2$ is the same as
the sum over square-free divisors of $a_3c_3$.  
But a simple application of partial summation now reveals that
\begin{align*}
\mcal{M}(A_i,B_i,C_i) 
&\ll \sum_{a_i,b_1,b_2,c_3} d(a_3)d(c_3)+ 
C_1C_2\sum_{a_i,b_1,b_2,c_3} \frac{d(a_3)d(c_3)}{|a_3c_3^2|}\\
&\ll \Big(AB_1B_2C_3+ A_1A_2B_1B_2C_1C_2\Big) (\log AC)^2,
\end{align*}
which thereby completes the proof of Proposition \ref{M2-x}.
\end{proof}

Although we shall not need to do so here, it is worth pointing out
that with more work it is possible to remove the term
$(\log AC)^2$ from the statement of Proposition \ref{M2-x}.

\section{Passage to the universal torsor}\lab{ut}

Our goal in this section is to equate the quantity (\ref{n(b)}) to the cardinality of a certain
subset of integral points on the universal torsor above $\tilde{X}$.
In fact our approach to the universal torsor rests upon an entirely elementary analysis of the
equation defining the surface $X$, and we shall not prove here that the
resulting parametrisation is actually the universal torsor above $\tilde{X}$.
This fact will be supplied for us by the work of Hassett and Tschinkel \cite{hassetttschinkel03}.

In any solution $\x \in Z_*^4$ to the equation $F(\x)=0$ we see that $x_4$
divides $x_1x_2x_3$.  Hence we may write $x_4=y_1y_2y_3$ and
$x_i=y_i z_i$, for some $\ma{y,z} \in \Z^3$ with $y_iz_i \neq 0$.   
Suppose that $z_i=\ve_iz_i'$ for $\ve_i =\pm 1$ and $z_i'\in \N.$  
Then one easily employs the equation $F(\x)=0$ to deduce that $\ve_1\ve_2\ve_3=1$.
Hence, upon relabelling variables we may assume that 
$$
x_i=y_i z_i, \quad x_4=y_1y_2y_3,
$$
for $(\ma{y,z})\in \Z^3\times\N^3$ with $y_i \neq 0$.  

Under this
substitution,    the equation $F(\x)=0$ becomes 
\beq\lab{eq1}
z_1z_2z_3=(y_1z_1+y_2z_2+y_3z_3)^2.
\eeq
Since $\x$ is assumed to be primitive, it follows that $\ma{y}$ is
primitive.  Moreover, if $p$ is any prime divisor of 
$\hcf(z_i,y_j)$, then (\ref{eq1}) implies that $p$ divides $y_kz_k$.
But this contradicts the primitivity of $\x$, whence 
\beq\lab{polo}
\hcf(z_i,y_j)=\hcf(y_1,y_2,y_3)=1.
\eeq
We now write
$z_i=w_it_i^2$, for square-free $w_i \in \N$ and non-zero $t_i \in
\Z$.  In fact we may assume that $t_i \in \N$, since $t_i$ and $-t_i$
produce the same value of $z_i$.

Then it follows from (\ref{polo}) that
\beq\lab{golf}
\hcf(w_i,y_j)=\hcf(t_i,y_j)=1,
\eeq
and from (\ref{eq1}) that $w_1w_2w_3$ is a square. Hence we can write
$$
w_1=u_2u_3, \quad w_2=u_1u_3, \quad w_3=u_1u_2,
$$
for square-free $u_i \in \N$, satisfying 
\beq\lab{eq2}
\hcf(u_i,u_j)=\hcf(u_i,y_i)=1.
\eeq 
Indeed $w_i$ is square-free, and any prime divisor of $\hcf(u_i,y_i)$
must also divide $\hcf(w_jw_k,y_i)$, contrary to (\ref{golf}).

Substituting the 
quantities $w_i=u_ju_k$ into (\ref{eq1}) therefore yields the expressions
\beq\lab{teq3}
\ve t_1t_2t_3u_1u_2u_3=y_1u_2u_3t_1^2 + y_2u_1u_3t_2^2 + y_3u_1u_2t_3^2,
\eeq
where $\ve=\pm 1$.
It is clear that $u_i$ must divide $y_iu_ju_kt_i^2$.  But then $u_i$ divides $t_i$,
since $u_i$ is square-free and $\hcf(u_i,y_iu_ju_k)=1$, by
(\ref{eq2}).  We proceed by writing
$$
s_0=\hcf(t_1/u_1,t_2/u_2,t_3/u_3),
$$
and $s_i=t_i/(s_0u_i)$.  Plainly $s_0, s_i \in \N$, and 
$\ma{s}=(s_1,s_2,s_3)$ is primitive.  Moreover, (\ref{golf}) yields 
\beq\lab{eq3}
\hcf(u_i,y_j)=\hcf(s_i,y_j)=1.
\eeq
Substituting $t_i=s_0s_iu_i$ into (\ref{teq3}), we therefore
obtain the equations
\beq\lab{ut1'}
\ve s_0s_1s_2s_3u_1u_2u_3 = y_1u_1s_1^2 + y_2u_2s_2^2 + y_3u_3s_3^2,
\eeq
where $\ve=\pm 1$, and $\hcf(s_0,y_i)=1$ by (\ref{golf}).

We proceed by using this equation, together with the fact that
$\ma{s}$ is primitive, to establish that 
$$
\hcf(s_i,u_j)=1.
$$
If $p$ is any prime divisor of $\hcf(s_i,u_j)$ then it follows
from (\ref{ut1'}), in conjunction with the coprimality conditions
(\ref{eq2}) and (\ref{eq3}),  that $p$ divides $s_k$. Considering
the corresponding $p$-adic order of each of the terms in (\ref{ut1'}),
one is easily led to the conclusion that $p$ divides $y_js_j^2$, since
$u_j$ is square-free.  But $\hcf(u_j,y_j)=1$ by (\ref{eq2}), and so 
$p$ divides $s_j$, which is impossible.
In fact we may go further and deduce that the components of $\ma{s}$
satisfy the relation
$$
\hcf(s_i,s_j)=1.
$$
This follows immediately from (\ref{eq3}), (\ref{ut1'}) and the fact
that $\hcf(s_i,u_k)=1$. 

Let 
$\mcal{T} \subset \A^{10}$ denote the set of $(s_0,\ma{s,u,y}) \in
\N \times N^3 \times N^3 \times Z_*^3$ such that 
\beq\lab{ut1}
s_0s_1s_2s_3u_1u_2u_3 = y_1u_1s_1^2 + y_2u_2s_2^2 + y_3u_3s_3^2,
\eeq
with 
\beq\lab{coprime1}
|\mu(u_1u_2u_3)|=\hcf(s_i,s_j)=\hcf(s_i,u_j)=1,
\eeq
and
\beq\lab{coprime2}
\hcf(s_0,y_i)=\hcf(s_i,y_j)=\hcf(u_i,y_1y_2y_3)=1.
\eeq
Now let $\x \in Z_*^4$ be any solution to the equation $F(\x)=0$.
Then tracing back through our argument, we deduce that there exists $(\pm s_0,\ma{s,u,y}) \in
\mcal{T}$ such that
\beq\lab{xi}
x_i=y_iu_i^2u_ju_ks_0^2s_i^2, \quad x_4=y_1y_2y_3.
\eeq
Conversely, given any $(\pm s_0,\ma{s,u,y}) \in \mcal{T}$, the point
given by 
(\ref{xi}) will be a primitive integer solution of the equation $F(\x)=0$, with
$x_1x_2x_3x_4 \neq 0$.  
Indeed if $p$ is any prime divisor of $x_1,x_2,x_3,x_4$ then we may assume
that 
$$
p\mid y_i, \quad p \mid s_0^2u_1u_2u_3\hcf(y_ju_js_j^2,y_ku_ks_k^2).
$$
But then (\ref{coprime2}) implies that $p \mid \hcf(y_1,y_2,y_3)$,
which is impossible.
We have therefore established the following result.

\begin{lem}\lab{base}
We have
$$
N_{U,H}(B)=\frac{1}{4}\#\Big\{(s_0,\ma{s,u,y})\in \mcal{T}:
  \max\{|y_iu_i^2u_ju_ks_0^2s_i^2|,
  |y_1y_2y_3|\} \leq B\Big\}.
$$
\end{lem}

The equation (\ref{ut1}) is an affine embedding of the universal torsor above the minimal
desingularisation $\tilde{X}$ of $X$. As already mentioned, it has been calculated by Hassett and
Tschinkel \cite[\S $4$]{hassetttschinkel03} by computing generators for the Cox
ring $\rom{Cox}(\tilde{X})$ of $\tilde{X}$.

\section{The lower bound}

Our method of establishing the lower bound closely follows 
Heath-Brown's treatment of the Cayley cubic.  
Consequently we shall adopt similar notation throughout this section.

Let $P_1,P_2 \leq B^\delta$, for some suitably small choice of
$\delta>0$.  This choice will be specified in (\ref{201}), below.
The idea is to fix choices of $\ma{s,u} \in N^3$ such that
(\ref{coprime1}) holds and 
$$
u_1u_2u_3=P_1, \quad s_1s_2s_3=P_2.
$$
In fact we shall insist upon the stronger condition that $P_1P_2$
is square-free.  This is clearly permissable for the purposes of a
lower bound.
We then  count the number of comparatively large non-zero solutions $s_0,y_1,y_2,y_3$ to
the linear equation (\ref{ut1})
subject to certain constraints.

Thus for  $Y_0,Y_i \geq 1$, we let 
$$
\mcal{N}=\mcal{N}(\ma{s,u};Y_0,Y_1,Y_2,Y_3)
$$ 
denote the number of $(s_0,y_1, y_2, y_3) \in \N\times \Z^3$ constrained by
(\ref{ut1}) and
\beq\lab{durham}
\hcf(s_0,y_i)= \hcf(y_i,P_1P_2)=1,
\eeq
for which 
\beq\lab{range-lower}
Y_0\leq s_0 <2Y_0, \quad Y_i\leq  |y_i| <2Y_i.
\eeq
It should be clear that whenever (\ref{ut1}) and (\ref{durham})  both hold, we
automatically have (\ref{coprime2}) and 
$\hcf(y_1,y_2,y_3)=1$.
It will be convenient to define the quantities
$$
A_0=P_1P_2, \quad
A_i=u_is_i^2,
$$
so that (\ref{ut1}) may be written
\beq\lab{P1}
A_0s_0=A_1y_1+A_2y_2+A_3y_3.
\eeq
Now it follows from Lemma \ref{base} that we are only interested in
values of $s_0,y_i$ for which 
$$
A_0^2A_i |s_0^2y_i| \leq BP_1P_2^2, \quad A_1A_2A_3|y_1y_2y_3| \leq BP_1P_2^2.
$$
Hence we shall choose 
\beq\lab{newcastle}
Y_0=\Big[\frac{(BP_1P_2^2)^{1/3}}{2A_0}\Big], \quad
Y_i=\Big[\frac{(BP_1P_2^2)^{1/3}}{2A_i}\Big]. 
\eeq

Much as in Heath-Brown's treatment, the main difficulty arises from
having to keep track of the coprimality conditions (\ref{durham}).  Let
$$
Q=P_1P_2 \prod_{p \leq \sqrt{\log B}}p.
$$
Following
\cite[Equation (3.6)]{cayley}, we  write
\beq\lab{diff}
\mcal{N}\geq \mcal{N}_1 - \mcal{N}_2,
\eeq
where $\mcal{N}_1$ is the number of solutions in which the condition
$\hcf(s_0,y_i)=1$ is replaced by the weaker condition
$$
\hcf(s_0,y_i,Q)=1,
$$
and $\mcal{N}_2$ is the number of solutions in which some $y_i$ shares
a prime factor $p$ with $s_0$, such that $p \nmid Q$.

We proceed by estimating $\mcal{N}_1$, for which we use the M\"obius
function to pick out the coprimality conditions.
Let $\mcal{N}_3(\ma{d};\ma{e})=\mcal{N}_3(d_1,d_2,d_3;e_{1},e_2,e_{3})$
denote the number of solutions of the equation (\ref{P1}) in the range
(\ref{range-lower}) with $d_i \mid y_i$ and  
$e_i \mid s_0,y_i$.  Then
\beq\lab{ben}
\mcal{N}_1=\sum_{d_i \mid P_1P_2}\mu(d_1)\mu(d_2)\mu(d_3)\sum_{e_{i} \mid
  Q}\mu(e_1)\mu(e_2)\mu(e_3) \mcal{N}_3(\ma{d};\ma{e}).
\eeq
Our task is to estimate $\mcal{N}_3(\ma{d};\ma{e})$. Define the least common
multiples
$$
h_0=[e_1,e_2,e_3], \quad h_i=[d_i,e_i].
$$
and the lattice
$$
\sfl=\{(n_1,n_2,n_3)\in \Z^3: A_ih_i \mid n_i, \;A_0h_0 \mid n_1+n_2+n_3\}.
$$
Then upon defining the region
$$
\mcal{R}=\{\ma{r} \in \R^3:  A_iY_i \leq |r_i|< 2A_iY_i, ~ A_0Y_0\leq |r_1+r_2+r_3|< 2 A_0Y_0\},
$$
one follows the lines of Heath-Brown's argument in order to deduce that
$$
\mcal{N}_3(\ma{d};\ma{e}) = \frac{\rom{vol}(\mcal{R})}{\det \sfl} +O((\det \sfl)^2 \max \{Y_i,Y_0\}^2),
$$
where
$$
\det \sfl = \frac{A_0h_0\prod_i A_ih_i}{\hcf(A_0h_0,A_ih_i)}.
$$
Since $d_i \mid P_1P_2$ and $e_i \mid Q$, we deduce that 
$A_0h_0 \leq P_1P_2Q^3$ and $A_ih_i
\leq P_1^2P_2^3Q$.  Hence we have 
$$
\det \sfl \leq P_1^7P_2^{10}Q^6 \ll P_1^{13}P_2^{16} \exp(O(\sqrt{\log
  B})) \ll B^{30\delta}.
$$ 
It follows that the
error term in our estimate for $\mcal{N}_3(\ma{d};\ma{e})$ is $O(B^{2/3+62\delta})$,
and so (\ref{ben}) becomes
\beq\lab{mp3}
\mcal{N}_1=\rom{vol}(\mcal{R}) 
\sum_{d_i, e_i}\mu(d_1)\cdots\mu(e_3) \frac{\hcf(A_0h_0,A_ih_i)}{A_0h_0\prod_i A_ih_i}
   + O(B^{2/3+63\delta}),
\eeq
since there are at most $O(B^\delta)$ divisors of $P_1P_2Q$.  

We now investigate the sum
\beq\lab{xmas}
\sum_{\colt{d_i \mid P_1P_2}{e_i\mid Q}}\mu(d_1)\cdots\mu(e_3)
\frac{\hcf(A_0h_0,A_ih_i)}{A_0h_0\prod_i A_ih_i}= \prod_{p\mid Q} E_p,
\eeq
say. When $p \mid Q$, but $p \nmid
P_1P_2$, we see that
$$
E_p=\sum_{\ve_i \geq 0} \mu(p^{\ve_1}) \mu(p^{\ve_2}) \mu(p^{\ve_3})
\frac{\hcf(p^{\max\{\ve_1,\ve_2,\ve_3\}},p^{\ve_i})}{
p^{\max\{\ve_1,\ve_2,\ve_3\}+\ve_1+\ve_2+\ve_3}},
$$
from which it easily follows that
\beq\lab{e1}
E_p=1-\frac{3}{p^2}+\frac{2}{p^3}.
\eeq
In this calculation we have used the fact that $p\nmid
A_0A_i$ whenever $p\nmid P_1P_2$.
Next, when  $p \mid P_1$ we may assume that $p$ divides precisely one factor, $u_1$
say. Since $P_1P_2$ is square-free it follows that $p \nmid
u_2u_3P_2$, and that $p^2 \nmid u_1$.  Let $A_0'=p^{-1}A_0, A_1'=p^{-1}A_1, A_2'=A_2$
 and $A_3'=A_3$, so that
$$
\frac{\hcf(A_0h_0,A_ih_i)}{A_0h_0\prod_i A_ih_i}=
\frac{\hcf(pA_0'h_0,pA_1'h_1,A_2'h_2,A_3'h_3)}{p^2A_0'h_0\prod_i
  A_i'h_i},
$$
with $p\nmid A_0'A_i'$.
Then in this setting we see that
$$
E_p=\frac{1}{p^2}\sum_{\delta_i, \ve_i \geq 0} \mu(p^{\delta_1})\mu(p^{\delta_2})\cdots \mu(p^{\ve_3})
\frac{\hcf(p,p^{\max\{\delta_2,\ve_2\}}, p^{\max\{\delta_3,\ve_3\}})}{
p^{\max\{\ve_1,\ve_2,\ve_3\}+\max\{\delta_1,\ve_1\}+\cdots+ \max\{\delta_3,\ve_3\}}},
$$
whence a straightforward calculation yields
\beq\lab{e2}
E_p=\frac{1}{p^2}\Big(1-\frac{1}{p}-\frac{1}{p^2}+\frac{1}{p^3}\Big).
\eeq
Finally we consider the case $p \mid P_2$, so that $p$ divides precisely one factor, $s_1$
say.  Since $P_1P_2$ is square-free it follows that $p \nmid
P_1s_2s_3$, and that $p^2 \nmid s_1$.  Let $A_0'=p^{-1}A_0,  A_1'=p^{-2}A_1, A_2'=A_2$
 and $A_3'=A_3$.  Then arguing as above we now have
$$
E_p=\frac{1}{p^3}\sum_{\delta_i, \ve_i \geq 0} \mu(p^{\delta_1})\mu(p^{\delta_2})\cdots \mu(p^{\ve_3})
\frac{\hcf(p,p^{\max\{\delta_2,\ve_2\}}, p^{\max\{\delta_3,\ve_3\}})}{
p^{\max\{\ve_1,\ve_2,\ve_3\}+\max\{\delta_1,\ve_1\}+\cdots+ \max\{\delta_3,\ve_3\}}}.
$$
In view of our calculation for (\ref{e2}) we immediately deduce that
\beq\lab{e3}
E_p=\frac{1}{p^3}\Big(1-\frac{1}{p}-\frac{1}{p^2}+\frac{1}{p^3}\Big).
\eeq
Taking (\ref{e1})-- (\ref{e3}) together in (\ref{xmas}), it therefore
follows that
$$
\sum_{d_i, e_i}\mu(d_1)\cdots\mu(e_3) \frac{\hcf(A_0h_0,A_ih_i)}{A_0h_0\prod_i A_ih_i} \gg
\frac{1}{P_1^2P_2^3}\frac{\phi(P_1P_2)}{P_1P_2}=\frac{\phi(P_1P_2)}{P_1^3P_2^4},
$$
since $\phi(n)=n\prod_{p\mid n}(1-1/p)$ for any $n \in \N$.
Our choices (\ref{newcastle}) for $Y_i,Y_0$ clearly imply that
$\rom{vol}(\mcal{R}) \gg BP_1P_2^2$.  We claim that
\beq\lab{lower-N1}
\mcal{N}_1 \gg \frac{B}{P_1P_2} \frac{\phi(P_1P_2)}{P_1P_2},
\eeq
provided that we take 
\beq\lab{201}
\delta=1/201.
\eeq
In order to establish the claim, it clearly suffices to check that the
lower bound in (\ref{lower-N1}) is larger than the error term in (\ref{mp3})
when $\delta$ is taken to  be $1/201$.
But on using the trivial lower bound $\phi(n)\geq 1$ for any $n
\in \N$, we see that
$$
\frac{B}{P_1P_2} \frac{\phi(P_1P_2)}{P_1P_2}
\geq \frac{B}{(P_1P_2)^2} \gg B^{1-4\delta}.
$$
Since $B^{1-4\delta}\gg B^{2/3+63\delta}$ for $\delta =1/201$, the
claim follows.

Next we must produce an upper bound for $\mcal{N}_2$, for which we may
ignore any coprimality conditions whenever we wish to.  Suppose that
$p\mid s_0, y_1$, for some prime $p \nmid Q$ lying in the range $R
\leq p <2R.$  In particular we may assume that $R \ll Y_1$.  There are
$O(R)$ such primes, and we fix one particular choice.
Following Heath-Brown's treatment, we write
$s_0=pt_0$ and $y_1=pt_1$ and count solutions of the linear equation
\beq\lab{p-linear}
pA_0t_0=pA_1t_1+A_2y_2+A_3y_3.
\eeq
In particular $t_0,t_1$ are contained in the ranges
$$
Y_0/R \ll |t_0| \ll Y_0/R, \quad Y_1/R \ll |t_1| \ll Y_1/R.
$$
Since $P_1P_2$ is square-free, it follows that $\hcf(A_i,A_0)=u_is_i$.
Hence we may deduce from (\ref{p-linear}) that 
$$
A_2y_2 \equiv -A_3y_3 \mod{pu_1s_1}.
$$
We may assume by symmetry that $A_2Y_2 \geq A_3Y_3$. Upon
noting that $pu_1s_1$ is coprime to $A_2$, since $\hcf(y_i,P_1P_2)=1$, it
follows that for each choice of $y_3$, there are
$O(1+Y_2/(Ru_1s_1))$ possibilities for $y_2$.  Now (\ref{p-linear})
implies that
\beq\lab{a1}
Ru_1s_1 \ll \max\{A_2Y_2,A_3Y_3\} = A_2Y_2.  
\eeq
Moreover, it follows from (\ref{newcastle}) that 
\beq\lab{a2}
A_2u_1s_1 \ll Y_2,
\eeq
provided that $\delta \leq 1/5$.  Together (\ref{a1}) and (\ref{a2})
imply that 
$$
1 \ll
\Big(\frac{A_2Y_2}{Ru_1s_1}\Big)^{1/2}\Big(\frac{Y_2}{A_2u_1s_1}\Big)^{1/2}=
\frac{Y_2}{R^{1/2}u_1s_1},
$$
whence we deduce that there are $O(Y_2Y_3/(R^{1/2}u_1s_1))$ choices for $y_2,y_3$.
We fix such a choice and write $A_2y_2+A_3y_3=pu_1s_1k$.  Then it
remains to count values of $t_0,t_1$ for which
\beq\lab{p-linear2}
u_2u_3s_2s_3t_0=s_1t_1+k.
\eeq
Now we have already seen that $R \ll Y_1$.  Moreover, as in
(\ref{a2}), we can use (\ref{newcastle}) to show that
$u_2^3u_3^3s_2^3s_3^3\ll Y_1$ provided that $\delta \leq 1/15$.
Together these inequalities imply that
$$
1 \ll
\Big(\frac{Y_1}{R}\Big)^{2/3}\Big(\frac{Y_1}{u_2^3u_3^3s_2^3s_3^3}\Big)^{1/3}=
\frac{Y_1}{R^{2/3}u_2u_3s_2s_3}.
$$
Viewing (\ref{p-linear2}) as a congruence modulo $u_2u_3s_2s_3$, one
easily concludes that there are $O(Y_1/(R^{2/3}u_2u_3s_2s_3))$
possibilities for $t_0,t_1$.  

In conclusion we have therefore shown that the total number of
admissible $p,y_2,y_3,t_0,t_1$, for which $R \leq p <2R$, is
$$
\ll R \cdot\frac{Y_2Y_3}{R^{1/2}u_1s_1}\cdot\frac{Y_1}{R^{2/3}u_2u_3s_2s_3} \ll
\frac{BP_1P_2^2}{R^{1/6}P_1^2P_2^3}= \frac{B}{R^{1/6}P_1P_2},
$$
by (\ref{newcastle}).  
Summing $R \gg \sqrt{\log B}$ over dyadic intervals, we deduce that
$$
\mcal{N}_2 \ll \frac{B}{P_1P_2}(\log B)^{-1/12},
$$
provided that $\delta \leq 1/15.$   It follows from (\ref{lower-N1})
and (\ref{201}) that $\mcal{N}_2=o(\mcal{N}_1)$, and so (\ref{diff}) implies that 
$$
\mcal{N}\gg \frac{B}{P_1P_2}\frac{\phi(P_1P_2)}{P_1P_2}.
$$

Finally, in order to complete the proof of the lower bound in
(\ref{theorem}), 
we note that any square-free value of $P$ will factorise
into values $u_1,u_2,u_3,s_1,s_2,s_3$ satisfying (\ref{coprime1}), in
precisely $d_6(P)$ ways.  It therefore follows that 
\beq\lab{train}
N_{U,H}(B) \gg \sum_{P\leq B^{2/201}}
|\mu(P)|d_6(P)\frac{B}{P}\frac{\phi(P)}{P}. 
\eeq
To handle this quantity we define the sum
$$
S(x)=\sum_{n \leq x} \frac{|\mu(n)|d_6(n)\phi(n)}{n},
$$
for any $x>1$, and proceed by establishing the following simple
bound.

\begin{lem}\lab{claim-xmas}
For any $x>1$ we have 
$$
S(x)\gg x(\log x)^5.
$$
\end{lem}

\begin{proof}
To establish the lemma we shall apply Perron's formula to the
corresponding Dirichlet series
$$
F(s)=\sum_{n=1}^\infty \frac{|\mu(n)|d_6(n)\phi(n)/n}{n^{s}},
$$
defined for $\Re e(s)>1$.  It is a trivial matter to calculate the
Euler product
$$
F(s)=\prod_p \Big(1+\frac{6(1-1/p)}{p^s}\Big)= \zeta(s)^6G(s),
$$
for some function $G(s)$ that is holomorphic and bounded on the
half-plane $\Re e(s)>1/2$.  Let $\ve>0$ and let $T \in [1,x]$.
Then Perron's formula yields
$$
S(x)= \frac{1}{2\pi i} \int_{1+\ve-iT}^{1+\ve+iT}\zeta(s)^6G(s) \frac{x^s}{s}\d s + O_\ve\Big( \frac{x^{1+\ve}}{T}\Big).
$$
We apply Cauchy's residue theorem to the rectangular contour 
joining the points ${2/3-iT}$, ${2/3+iT}$,
${1+\ve+iT}$ and ${1+\ve-iT}$, which therefore leads to the conclusion that
there exists a polynomial $f$ of degree $5$ such that
$$
S(x)-xf(\log x) \ll_\ve 
\frac{x^{1+\ve}}{T}+
\Big(\int_{2/3-iT}^{2/3+iT}+\int_{2/3-iT}^{1+\ve-iT}+ 
\int_{1+\ve+iT}^{2/3+iT}\Big) \Big|\zeta(s)^6\frac{x^s}{s}\Big|\d s.
$$
Here we have used the fact that $G(s)$ is bounded for $\Re e(s)>1/2.$
To estimate this error term we apply the well known convexity bound
$\zeta(\sigma+i t)\ll_\ve |t|^{(1-\sigma)/3+\varepsilon},$
valid for any $\sigma\in[1/2,1]$ and $|t|\geq 1$.  Thus it follows that
$$
S(x)-xf(\log x) \ll_\ve 
\frac{x^{1+\ve}}{T}+ x^{2/3+\ve}T^{2/3}.
$$
Selecting $T=x^{1/5}$ therefore completes the proof of Lemma \ref{claim-xmas}.
\end{proof}

On combining Lemma \ref{claim-xmas} with an application of partial
summation, and then inserting the resulting estimate into
(\ref{train}), we therefore deduce that 
$$
N_{U,H}(B) \gg B (\log B)^6.
$$
This completes the proof of the lower bound in (\ref{theorem}).

\section{The upper bound}\lab{s-upper}

Fix a choice of $X_1,\ldots,X_4, S_0,S_i, U_i, Y_i \geq 1$. 
We shall write 
$$
\mcal{N}=\mcal{N}(X_1,\ldots,X_4;S_0;S_1,S_2,S_3;U_1,U_2,U_3;Y_1,Y_2,Y_3)
$$
for the total contribution to $N_{U,H}(B)$ from $\x$ contained in the intervals
\beq\lab{range-xi}
X_\xi\leq |x_\xi| < 2X_\xi, \quad (1\leq \xi \leq 4),
\eeq
and such that the variables $s_0,\ma{s,u,y}$ 
appearing in Lemma \ref{base} satisfy
\beq\lab{range1}
S_0\leq s_0 < 2S_0, \quad S_i\leq |s_i| < 2S_i, \quad U_i\leq |u_i| <2U_i, \quad
Y_i\leq |y_i| <2Y_i.
\eeq
It will be convenient to relable the indices so that
\beq\lab{order}
X_1\leq X_2 \leq X_3.
\eeq
Suppose that $\x \in Z_*^4$ is
a solution of $F(\x)=0$, with $|x_1|,\ldots,|x_4| \leq B$.
Then (\ref{range-xi}) implies that
\beq\lab{Xi}
X_1, X_2, X_3, X_4 \leq B.
\eeq
If $\mcal{N}=0$ there is nothing to prove, and so we assume henceforth
that the dyadic ranges in (\ref{range-xi}) and (\ref{range1}) produce
a non-zero value of $\mcal{N}$.

We proceed by showing that under the assumption that $\mcal{N}\neq 0$,
certain choices of dyadic ranges in
(\ref{range-xi}) and (\ref{range1}) force certain other ranges to
have fixed order of magnitude.  
It will be convenient to write
$$
S=S_1S_2S_3, \quad U=U_1U_2U_3, \quad Y=Y_1Y_2Y_3.
$$
Hence it follows from (\ref{xi}) that 
\beq\lab{Yi}
X_i\ll Y_iU_iUS_0^2S_i^2 \ll X_i,
\eeq
and that
\beq\lab{Y}
X_4 \ll Y \ll X_4.
\eeq
Together, (\ref{Yi}) and (\ref{Y})  imply that
\beq\lab{S0}
\Big(\frac{X_1X_2X_3}{X_4}\Big)^{1/2} \ll S_0^3SU^2\ll \Big(\frac{X_1X_2X_3}{X_4}\Big)^{1/2}.
\eeq
We take a moment to record two further inequalities satisfied by the
quantities $S_0,S_i,U_i,Y_i$, which will be crucial in our final
analysis. First we deduce from 
(\ref{Yi}), (\ref{Y}) and (\ref{S0}) that
\begin{align}
S_0S^{1/3}U^{2/3}Y^{2/3} &= (S_0^3SU^{2})^{1/3}Y^{2/3}\nonumber\\
&\ll (X_1X_2X_3)^{1/6}X_4^{1/2}.\lab{pisa}
\end{align}
Similarly, we may deduce that
\begin{align}
S_0S^{1/2}UY^{1/2} &\leq
(S_0^{3}SU^2)^{1/2}Y^{1/2}\nonumber\\
&\ll (X_1X_2X_3X_4)^{1/4}.
\lab{florence}
\end{align}

It is clear that $\mcal{N}$ is bounded above by the number
of $s_0 \in \N$ and $\ma{s}, \ma{u}, \ma{y} \in Z_*^3$ 
contained in the ranges (\ref{range1}),  for
which (\ref{ut1}), (\ref{coprime1}) and (\ref{coprime2})  all hold.
Ultimately we shall sum over suitable dyadic intervals for $X_1,X_2,X_3,X_4$ 
and $S_0,S_i, U_i,Y_i$, in order to establish the 
upper bound in (\ref{theorem}), and so complete the proof of the theorem.

For any fixed choice of $s_0 \in \N$ and $\ma{s}, \ma{u} \in Z_*^3$ in the
region (\ref{range1}), with (\ref{coprime1}) holding, we let 
$$
N(s_0,\ma{s,u})=N(Y_1,Y_2,Y_3;s_0,\ma{s,u})
$$ 
denote the corresponding contribution to $\mcal{N}$
from the $\ma{y} \in Z_*^3$.
Clearly we are only interested in values of $s_0,\ma{s,u}$ for which
$N(s_0,\ma{s,u})$ is non-zero.  
Considering $s_0, \ma{s,u}$ to be fixed, we select any vector
\beq\lab{bartok}
\hat{\ma{y}}=(\hat{y_1},\hat{y_2},\hat{y_3})
\eeq
for which the Euclidean norm $|\hat{\ma{y}}|$ is least.  Following
the convention that this vector too is fixed, for fixed values of
$s_0,\ma{s},\ma{u}$, we define the change of variables
\beq\lab{zi}
z_i=y_i-\hat{y_i}.
\eeq
We shall let $N_1(s_0,\ma{s,u};\hat{\ma{y}})$ denote the overall contribution to
$N(s_0,\ma{s,u})$ from those $\ma{y}$ for which $z_1z_2z_3 \neq 0$,
and we let $N_2(s_0,\ma{s,u};\hat{\ma{y}})$ denote the remaining contribution to
$N(s_0,\ma{s,u})$ from those $\ma{y}$ for which $z_1z_2z_3 =0$.
With this notation we therefore have
\begin{align}
\mcal{N}
= \sum_{s_0,\ma{s,u}} N(s_0,\ma{s,u}) 
&= \sum_{s_0,\ma{s,u}} N_1(s_0,\ma{s,u};\hat{\ma{y}}) +\sum_{s_0,\ma{s,u}} N_2(s_0,\ma{s,u};\hat{\ma{y}})\nonumber\\
&=\mcal{N}_1+\mcal{N}_2,\lab{N}
\end{align}
say.  Here the summations are over all 
$(s_0, \ma{s,u}) \in \N\times Z_*^3\times Z_*^3$ in the region (\ref{range1}), with (\ref{coprime1}) holding.
It will be necessary to investigate the quantities $\mcal{N}_1$ and $\mcal{N}_2$ separately.
Finally we shall conclude that
\beq\lab{N'}
N_{U,H}(B)\ll \sum_{{S_0,S_i,U_i,Y_i,X_i, X_4}} \mcal{N},
\eeq
where the summation is over all dyadic intervals for 
$X_1,X_2,X_3,X_4$,  subject to 
(\ref{Xi}), and also all dyadic intervals for $S_0,S_i,U_i,Y_i$,
subject to (\ref{Yi}) and (\ref{Y}).
Our first task is to show that the overall contribution from
$\mcal{N}_1$ to $N_{U,H}(B)$ is satisfactory.

\begin{pro}\lab{n1}
We have
$$
\sum_{{S_0,S_i,U_i,Y_i,X_i, X_4}}
\mcal{N}_1 \ll B (\log B)^6,
$$
where the summation is over
dyadic intervals subject to (\ref{Xi})--(\ref{Y}).
\end{pro}

Proposition \ref{n1} will be established in \S \ref{char}.
Next in \S \ref{char2} we shall estimate the corresponding contribution from
$\mcal{N}_2$ to $N_{U,H}(B)$ via the following result.

\begin{pro}\lab{n2}
We have
$$
\sum_{{S_0,S_i,U_i,Y_i,X_i, X_4}}
\mcal{N}_2 \ll B (\log B)^6,
$$
where the summation is over
dyadic intervals subject to 
(\ref{Xi})--(\ref{Y}).
\end{pro}

Once taken together in (\ref{N}) and (\ref{N'}), Propositions \ref{n1} and \ref{n2}
therefore yield the upper bound
$$
N_{U,H}(B)\ll B(\log B)^6,
$$
which thereby completes the proof of (\ref{theorem}).

\subsection{Proof of Proposition \ref{n1}}\lab{char}

Our first step in the proof of Proposition \ref{n1} is to provide a pair of
upper bounds for $\mcal{N}_1$.
For any fixed choice of $s_0 \in \N$ and $\ma{s}, \ma{u} \in Z_*^3$ in the
region (\ref{range1}), with (\ref{coprime1}) holding, we let
$\hat{\ma{y}}$ be the corresponding vector (\ref{bartok}) that was
selected above.  On recalling the change of variables (\ref{zi}), it
therefore follows from (\ref{ut1})  that
\beq\lab{G}
z_1 u_1s_1^2+z_2 u_2s_2^2+ z_3 u_3s_3^2=0,
\eeq
and from (\ref{range1}) that $|z_i| <4Y_i$.
Hence we deduce that
$$
\mcal{N}_1
\ll 
S_0\sum_{\ma{s,u}}\#\{\ma{z} \in \Z^3: z_1z_2z_3\neq 0, ~|z_i|
<4Y_i, ~\mbox{(\ref{G}) holds}\}.
$$
Recall  that $u_1u_2u_3$ is square-free, so that $\hcf(u_i,u_j)=1$.
It is apparent that the $\ma{z}$ appearing in the summand need not be
primitive.  Moreover we no longer necessarily have coprimality
conditions corresponding to (\ref{coprime2}).    In order to recover a
weaker set of coprimality relations, we shall write
$$
u_i=d_{jk}u_i', \quad z_i=d_{ij}d_{ik}ez_i',
$$
say, for any  $d_{ij}, e \in \N$ with the convention that $d_{ij}=d_{ji}$. 
Let 
$$
U_i'=\frac{U_i}{d_{jk}}, \quad Y_i'=\frac{Y_i}{d_{ij}d_{ik}e},
$$ 
and 
$$
U'=U_1'U_2'U_3', \quad Y'=Y_1'Y_2'Y_3'. 
$$
In particular it follows from (\ref{coprime1}) that 
\beq\lab{cut}
|\mu(u_1'u_2'u_3')|=\hcf(s_i,s_j)=\hcf(s_i,u_j')=1,
\eeq
for any $\ma{s}$ and $\ma{u'}$.

For fixed values of $d_{ij}, e \in \N$, our task is to
estimate the number of $\ma{s,u',z'} \in Z_*^3$ such that 
$$
\hcf(u_i',z_j',z_k')=1,
$$
(\ref{cut}) holds, 
$$
S_i \leq |s_i| < 2S_i, \quad U_i' \leq |u_i'| < 2U_i',  \quad |z_i'| <4Y_i', 
$$
and
$$
z_1' u_1's_1^2+z_2' u_2's_2^2+ z_3'u_3's_3^2=0.
$$
But this quantity is clearly bounded above by
$\mcal{M}=\mcal{M}(2U_i',4Y_i',2S_i)$ in the notation of \S \ref{core-x}.
Thus it follows that
\beq\lab{ken}
\mcal{N}_1 \ll S_0\sum_{d_{ij}, e \ll Y}\mcal{M},
\eeq
and Proposition \ref{M1-x} yields
\beq\lab{M1}
\mcal{M} \ll_\ve S^{1/3}U'^{2/3}Y'^{2/3} +
\sigma\tau S^{1/2}U'Y'^{1/2},
\eeq
for any $\ve>0$, where
$$
\sigma=1+\frac{\min\{U, Y\}^\ve}{\min\{Y_i'Y_j'\}^{1/16}},
\quad 
\tau=1+\frac{\log Y}{\min\{Y_i'Y_j'\}^{1/16}}.
$$
On applying Proposition \ref{M2-x} we obtain the alternative estimate
\begin{equation}\lab{M2}
\mcal{M} \ll U'Y_i'Y_j'\Big(S_k +  S_iS_jU_k'^{-1}\Big) (\log SU')^2,
\end{equation}
for any permutation $\{i,j,k\}$ of the set $\{1,2,3\}$.

We may now use (\ref{M1}) and (\ref{M2}) to obtain a pair of  estimates for $\mcal{N}_1$.  Recall the
inequality (\ref{ken}) for $\mcal{N}_1$, and note that
$$
U'=(d_{12}d_{13}d_{23})^{-1}U, \quad
Y'=(d_{12}d_{13}d_{23})^{-2}e^{-3}Y. 
$$  
Beginning with an application of (\ref{M1}), we deduce that 
$$
\sum_{d_{ij}, e \ll Y}\mcal{M}\ll_\ve
S^{1/3}U^{2/3}Y^{2/3}
+\sigma\tau S^{1/2}UY^{1/2},
$$
where 
\beq\lab{sigtau}
\sigma=1+
\frac{\min\{U, Y\}^\ve}{\min \{Y_iY_j\}^{1/16}},
\quad
\tau=1+
\frac{\log B }{\min \{Y_iY_j\}^{1/16}}.
\eeq
We therefore obtain the following estimate.

\begin{lem}\lab{n1-1-lemma}
We have
$$
\mcal{N}_1 \ll_\ve  S_0S^{1/3}U^{2/3}Y^{2/3}
+\sigma\tau S_0S^{1/2}UY^{1/2},
$$
for any $\ve>0$, where $\sigma, \tau$ are given by (\ref{sigtau}).
\end{lem}

Similarly, since (\ref{Xi}) and (\ref{Yi}) imply that $(\log SU) \ll\log B$, 
an application of (\ref{M2}) yields the following alternative estimate.

\begin{lem}\lab{n1-2-lemma}
We have
$$
\mcal{N}_1 \ll
S_0UY_iY_j(S_k +  S_iS_jS_k^{-1}U_k^{-1}) (\log B)^2,
$$
for any permutation $\{i,j,k\}$ of the set $\{1,2,3\}$.
\end{lem}

We are now ready to complete the proof of Proposition \ref{n1}.
Now it follows from the
inequalities (\ref{range1}), (\ref{order}) and (\ref{Yi}) that
\beq\lab{1<2<3}
Y_1U_1S_1^2 \ll Y_2U_2S_2^2 \ll Y_3U_3S_3^2.
\eeq
In particular (\ref{ut1}) implies that
\beq\lab{beer}
S_0SU \ll Y_3U_3S_3^2.
\eeq
Multiplying both sides of this inequality by $S_0^2U$, and recalling (\ref{Xi}) and
(\ref{Yi}), we deduce that
\beq\lab{doctor}
S_0^3SU^2 \ll B.
\eeq
It will also be useful to deduce an inequality   involving the maximum
size of the $U_i$.   Suppose temporarily that
$U_1 \leq U_2 \leq U_3$, so that $U_1U_2 \leq U_3^2$.
Then in view of (\ref{Xi}) and (\ref{Yi}) it clearly
follows that $U_1U_2U_3^2 \ll B$, whence $U_1U_2 \ll B^{1/2}$.   Using this
sort of argument it is not hard to deduce that in general
\beq\lab{u-upper}
U_iU_j \ll B^{1/2}.
\eeq
Throughout the proof of Proposition \ref{n1} we shall
make the additional assumption that 
$$
Y_{\imath} \leq Y_{\jmath} \leq Y_{\kappa},
$$
for some permutation $\{\imath,\jmath,\kappa\}=\{1,2,3\}$.  Our plan is to
use Lemma \ref{n1-1-lemma} whenever 
\beq\lab{trout1}
Y_{\kappa} \leq  (Y_\imath Y_\jmath)^{9} \quad \mbox{or} \quad U \leq
(Y_\imath Y_\jmath)^{3},
\eeq
and Lemma \ref{n1-2-lemma} in the alternative case 
\beq\lab{trout2}
 Y_\imath Y_\jmath \leq  \min\{Y_\kappa^{1/9},U^{1/3}\}.
\eeq

Let us consider the case (\ref{trout1}) first.  Recalling the
definition (\ref{sigtau}) of $\sigma$, it follows that 
$$
\sigma \ll 1+ (Y_\imath Y_\jmath)^{10\ve-1/16}.
$$
Hence we may take
$\sigma \ll 1$ in Lemma \ref{n1-1-lemma}, provided that we take $\ve>0$
to be sufficiently small.
It therefore follows from (\ref{pisa}) and Lemma \ref{n1-1-lemma} that
\beq\lab{n1-final}
\mcal{N}_1\ll (X_1X_2X_3)^{1/6}X_4^{1/2} + 
S_0S^{1/2}UY^{1/2}\Big(1+ \frac{\log B}{(Y_\imath Y_\jmath)^{1/16}}\Big),
\eeq
whenever (\ref{trout1}) holds.
We now sum over the various dyadic intervals for 
${S_0,S_i,U_i,Y_i,X_i, X_4}$  subject to 
(\ref{Xi})--(\ref{Y}) and (\ref{trout1}).
Suppose for the moment that we want to sum over all
possible dyadic intervals $X \leq |x|<2X$, for which $|x| \leq
\mcal{X}$.  Then there are plainly $O(\log \mcal{X})$ possible choices
for $X$.  In addition to this basic estimate, we shall make frequent
use of the estimates 
$$
\sum_{X} X^{\delta} \ll_\delta 
\left\{
\begin{array}{ll}
1, & \delta<0,\\
\mcal{X}^{\delta}, & \delta>0.
\end{array}
\right.
$$

Returning to (\ref{n1-final}), we may
deduce from (\ref{Yi}) and (\ref{S0}) that values
of $S_0,Y_1,Y_2,Y_3$ are determined by the choices of 
$X_1,X_2,X_3,X_4$ and $S_i,U_i$.  
Now there  are clearly $\ll (\log B)^6$
possible sets of values  for 
$S_i,U_i$.  In view of (\ref{Xi}), we therefore obtain the estimate
\beq\lab{final-1}
\sum_{S_0,S_i,U_i,Y_i,X_i,X_4} (X_1X_2X_3)^{1/6}X_4^{1/2} \ll B (\log B)^6.
\eeq
Employing (\ref{florence}), we find similarly that 
\begin{align}\nonumber
\sum_{S_0,S_i,U_i,Y_i,X_i,X_4} S_0S^{1/2}UY^{1/2} &\ll
\sum_{S_0,S_i,U_i,Y_i,X_i,X_4} (X_1X_2X_3X_4)^{1/4}\\ 
&\ll B (\log B)^6.\lab{final-1'}
\end{align}
Finally we turn to the term 
$S_0S^{1/2}U(Y_\imath Y_\jmath)^{7/16}Y_\kappa^{1/2}\log B$ in (\ref{n1-final}).
We shall sum over dyadic intervals subject to the two inequalities
$$
Y_\kappa \ll \frac{B}{Y_\imath Y_\jmath}, \quad S_1 \ll \frac{B}{S_0^3S_2S_3U^2}.
$$
The first of these follows from (\ref{Xi}) and (\ref{Y}), whereas the
second is just (\ref{doctor}).   We therefore obtain the estimate
\begin{align*}
\sum_{S_0,S_i,U_i,Y_i} S_0S^{1/2}U(Y_\imath Y_\jmath)^{7/16}Y_\kappa^{1/2}
&\ll B^{1/2}\sum_{S_0,S_i,U_i,Y_\imath, Y_\jmath}
\frac{S_0S^{1/2}U}{(Y_\imath Y_\jmath)^{1/16}}\\
&\ll B\sum_{S_0,S_2,S_3,U_i,Y_\imath, Y_\jmath}
S_0^{-1/2}(Y_\imath Y_\jmath)^{-1/16}\\
&\ll B (\log B)^5.
\end{align*}
Since values of $X_1,X_2,X_3,X_4$ are 
determined by choices of  $S_0,S_i,U_i,Y_i$, we may 
combine this latter estimate with (\ref{final-1}) and  (\ref{final-1'})
in (\ref{n1-final}), in order to  conclude that 
\beq\lab{sum-N1-final}
\sum_{\colt{S_0,S_i,U_i,Y_i,X_i, X_4}{\mbox{\scriptsize{(\ref{trout1}) holds}}}}
\mcal{N}_1 \ll B (\log B)^6.
\eeq

Next we handle the case in which (\ref{trout2}) holds.
For this we employ the alternative estimate Lemma \ref{n1-2-lemma} to deduce that 
$$
\mcal{N}_1\ll (S_0S_\kappa UY_\imath Y_\jmath  + 
S_0S_\imath S_\jmath U_\imath U_\jmath Y_\imath Y_\jmath  )(\log B)^2.
$$
Again we proceed by summing over dyadic intervals for 
$S_0,S_i,U_i,Y_i,X_i,X_4$, this time subject to 
(\ref{Xi})--(\ref{Y}) and (\ref{trout2}).
Let us consider the term
$S_0S_\kappa U Y_\imath Y_\jmath $.  But then (\ref{Yi}), (\ref{S0}), (\ref{u-upper}) and
(\ref{trout2})  together imply that
\begin{align*}
S_0S_\kappa U Y_\imath Y_\jmath  &=  S_0^3SU^2 \frac{Y_\imath Y_\jmath }{S_0^2S_\imath S_\jmath U}\\ 
&\ll
(X_1X_2X_3)^{1/2}
\frac{(Y_\imath Y_\jmath )^{3/2}(U_\imath U_\jmath )^{1/2}}{(X_\imath X_\jmath )^{1/2}}\\ 
&\ll X_\kappa^{1/2} Y_\kappa^{1/6}  (U_\imath U_\jmath )^{1/2}\\
&\ll B^{11/12}.
\end{align*}
Since there are at most $O_\ve(B^\ve)$ 
dyadic intervals for $S_0,S_i,U_i,Y_i$, 
which in turn determine values of $X_1,X_2,X_3,X_4$,
this therefore  
leads to the conclusion that 
\beq\lab{toast} 
\sum_{{S_0,S_i,U_i,Y_i,X_i,X_4}}
S_0S_\kappa UY_\imath Y_\jmath  (\log B)^2 \ll B,
\eeq
whenever (\ref{trout2}) holds.
Lastly we consider the term $S_0S_\imath S_\jmath U_\imath U_\jmath Y_\imath Y_\jmath $.
Now there are  $O(\log B)$ dyadic intervals for $Y_\kappa$, and
(\ref{trout2}) implies that $Y_\imath ,Y_\jmath  \leq U^{1/3}$.
Employing the upper bound $S_\imath  \ll B/(S_0^3S_\jmath S_\kappa U^2)$, we 
therefore deduce that
\begin{align*}
\sum_{S_0,S_i,U_i,Y_i,X_i,X_4}
S_0S_\imath S_\jmath U_\imath U_\jmath Y_\imath Y_\jmath
&\ll \log B 
\sum_{S_0,S_i,U_i}
S_0S_\imath S_\jmath U_\imath U_\jmath U^{2/3}\\
&\ll 
B \log B 
\sum_{S_0,S_\jmath, S_\kappa ,U_i} S_0^{-2}S_\kappa^{-1}U^{-1/3}\\
&\ll
B (\log B)^2,
\end{align*}
whenever (\ref{trout2}) holds.
Once combined with (\ref{toast}) this yields the overall contribution
\beq\lab{sum-N1-final'}
\sum_{\colt{S_0,S_i,U_i,Y_i,X_i, X_4}{\mbox{\scriptsize{(\ref{trout2}) holds}}}}
\mcal{N}_1 \ll B (\log B)^4.
\eeq

Once taken together, (\ref{sum-N1-final}) and (\ref{sum-N1-final'})
therefore complete the proof of Proposition \ref{n1}.

\subsection{Proof of Proposition \ref{n2}}\lab{char2}

We begin this section by providing an upper bound for $\mcal{N}_2$.
For any fixed choice of $s_0 \in \N$ and $\ma{s}, \ma{u} \in Z_*^3$ in the
region (\ref{range1}), with (\ref{coprime1}) holding, let
$\hat{\ma{y}}$ be the vector (\ref{bartok}) counted by 
$N(s_0,\ma{s,u})$ that was selected at the start of \S \ref{s-upper}.   
Then (\ref{N}) implies that
$$
\mcal{N}_2
= \sum_{s_0,\ma{s,u}} N_2(s_0,\ma{s,u};\hat{\ma{y}}),
$$
where $N_2(s_0,\ma{s,u};\hat{\ma{y}})$ denotes the contribution to
$N(s_0,\ma{s,u})$ from those $\ma{y}$ for which 
$$
\prod_{1\leq i \leq
3}(y_i-\hat{y_i}) =0.
$$

Let $N_2^{(i)}$ denote the total contribution to $N_2(s_0,\ma{s,u};\hat{\ma{y}})$ from those 
$\ma{y}$ for which $y_i=\hat{y_i}$ is fixed.  
It therefore follows that
\beq\lab{sunday}
\mcal{N}_2 \leq \sum_{s_0,\ma{s,u}} \Big(N_2^{(1)}+N_2^{(2)}+N_2^{(3)}\Big)=
\mcal{N}_2^{(1)}+\mcal{N}_2^{(2)}+\mcal{N}_2^{(3)},
\eeq
say. In order to estimate $N_2^{(i)}$ for fixed values of $s_0\in \N$ and 
$\ma{s,u} \in Z_*^3$, it suffices to count non-zero integer solutions
$y_j,y_k$ to the equation
\beq\lab{n2-main}
y_ju_js_j^2+y_ku_ks_k^2=n, 
\eeq
where $n=s_0s_1s_2s_3u_1u_2u_3- \hat{y_i}u_is_i^2$ is fixed.  
Our first step is to deduce from (\ref{coprime1}) that
$$
\hcf(u_is_i^2,u_js_j^2)=1.
$$
Noting that $|y_j|<
2Y_j$ and $|y_k|< 2Y_k$, we proceed by applying  Lemma \ref{line} to (\ref{n2-main}).
Taking 
$$
\ma{h}=(u_js_j^2,u_ks_k^2,n), \quad \ma{w}=(y_j,y_k,1),
$$
we therefore deduce that 
\beq\lab{app-line}
N_2^{(i)}\ll 1+ \frac{Y_jY_k}{\max\{Y_jU_jS_j^2, Y_kU_kS_k^2,|n|\}}.
\eeq
Since $Y_iU_iS_i^2 \ll \hat{y_i}u_is_i^2 \ll Y_iU_iS_i^2$ and $S_0SU
\ll s_0s_1s_2s_3u_1u_2u_3 \ll S_0SU$, by  (\ref{range1}), it is easy
to see that
\begin{align*}
|n|&=|s_0s_1s_2s_3u_1u_2u_3-\hat{y_i}u_is_i^2|\\
&\geq \Big||\hat{y_i}u_is_i^2|-|s_0s_1s_2s_3u_1u_2u_3|\Big|\\
&\gg Y_iU_iS_i^2,
\end{align*}
if $Y_iU_iS_i^2 \gg S_0SU$.
Upon summing (\ref{app-line}) over all $s_0,\ma{s,u}$, and then
inserting the resulting bound into (\ref{sunday}), we therefore obtain
the following result

\begin{lem}\lab{for-n2-lemma}
We have
$$
\mcal{N}_2 \ll S_0SU + \max_{\{i,j,k\}}\Big\{\frac{S_0SU Y_jY_k}{\max\{Y_jU_jS_j^2,
  Y_kU_kS_k^2, \theta_i\}}\Big\},
$$
where the first maximum is over all permutations $\{i,j,k\}$ of the
set $\{1,2,3\}$, and 
\beq\label{theta}  
\theta_i= 
\left\{
\begin{array}{ll}
Y_iU_iS_i^2, & Y_iU_iS_i^2 \gg S_0SU,\\
1, & \mbox{otherwise.}
\end{array}
\right.
\eeq
\end{lem}

We now complete the proof of Proposition \ref{n2}.
Recall the inequality (\ref{beer}).
Our first task will be to establish that
\beq\lab{final-2}
\sum_{S_0,S_i,U_i,Y_i,X_i, X_4}  Y_3U_3S_3^2 \ll B (\log B)^4,
\eeq
where the summation is over dyadic intervals subject to  (\ref{Xi})--(\ref{Y}).
In order to do so we observe as in \S \ref{char} that values
of $X_1,X_2,X_3,X_4$ are 
determined by the choices of 
$S_0,S_i,U_i,Y_i$.
Recall the inequalities (\ref{1<2<3}).
We have two basic cases to consider, according to whether or not $Y_3U_3S_3^2$
is sufficiently large compared with $Y_2U_2S_2^2$.

Suppose first that $Y_3U_3S_3^2 \gg Y_2U_2S_2^2$. Then the ranges (\ref{range1})
imply that
$$
|y_1u_1s_1^2+y_2u_2s_2^2+y_3u_3s_3^2| \geq 
||y_3u_3s_3^2|-|y_1u_1s_1^2+y_2u_2s_2^2|| \gg Y_3U_3S_3^2,
$$ 
in any solution.
Since we obviously have $|y_1u_1s_1^2+y_2u_2s_2^2+y_3u_3s_3^2| \ll
Y_3U_3S_3^2$, the basic equation (\ref{ut1}) implies that
$Y_3U_3S_3^2 \ll S_0SU \ll Y_3U_3S_3^2,$
whence
\beq\lab{y3s3}
S_0S_1S_2U_1U_2 \ll Y_3S_3 \ll S_0S_1S_2U_1U_2.
\eeq
Summing over $Y_3 \ll S_0S_1S_2S_3^{-1}U_1U_2$, we therefore obtain
$$
\sum_{S_0,S_i,U_i,Y_i,X_i, X_4}  Y_3U_3S_3^2 
\ll \sum_{S_0,S_i,U_i,Y_1,Y_2}  S_0SU,
$$
where the last sum is subject to the inequality (\ref{doctor}).
Since there are $\ll (\log B)^2$ choices for $Y_1,Y_2$, we therefore see
that this sum is at most
\begin{align*}
&\ll (\log B)^2  
\sum_{S_0,S_1,S_2,U_i}  S_0S_1S_2U \sum_{S_3 \ll B/(S_0^3S_1S_2U^2)} S_3\\
&\ll B(\log B)^2  
\sum_{S_0,S_1,S_2,U_i}  S_0^{-2}U_1^{-1}U_2^{-1}U_3^{-1}\\
&\ll B (\log B)^4,
\end{align*} 
as required for (\ref{final-2}).

Next, if $Y_2U_2S_2^2 \ll Y_3U_3S_3^2 \ll Y_2U_2S_2^2$, then it
follows that any choice of $Y_3,U_2,U_3,S_2,S_3$ determines a choice
of $Y_2$.  Proceeding in a similar fashion to above, we deduce from
(\ref{Xi}) and  (\ref{Yi}) that $Y_3 \ll B/(U_3US_0^2S_3^2)$.  Hence
we obtain the estimate
\begin{align*}
\sum_{S_0,S_i,U_i,Y_i,X_i, X_4}  Y_3U_3S_3^2 & \ll B
\sum_{S_0,S_i,U_i,Y_1} S_0^{-2}U_1^{-1}U_2^{-1}U_3^{-1}\\
&\ll B (\log B)^4.
\end{align*}
This completes the proof of (\ref{final-2}).

Recall the estimate in Lemma \ref{for-n2-lemma} for $\mcal{N}_2$. Then
in view of (\ref{final-2}),  it suffices to estimate
\beq\lab{sf}
N^{(i)}(B)=\sum_{S_0,S_i,U_i,Y_i} \frac{S_0SU Y_jY_k}{\max\{Y_jU_jS_j^2,
  Y_kU_kS_k^2, \theta_i\}},
\eeq
for each permutation $\{i,j,k\}$ of $\{1,2,3\}$, and where $\theta_i$ is given by (\ref{theta}).
We begin by handling the case $i=3$. 
Suppose first that $Y_3U_3S_3^2 \gg Y_2U_2S_2^2$, so that (\ref{y3s3})
holds and we may take $\theta_3=Y_3U_3S_3^2$.  Then 
$$
\frac{S_0SU Y_1Y_2}{\max\{Y_1U_1S_1^2,
  Y_2U_2S_2^2, \theta_3\}} = \frac{S_0S_1S_2U_1U_2Y_1Y_2}{S_3Y_3}.
$$
Moreover, we recall the inequalities
$$
Y_2\ll \frac{B}{U_2US_0^2S_2^2}, \quad Y_1 \ll
\frac{Y_3U_3S_3^2}{U_1S_1^2}, \quad S_3 \ll \frac{S_0S_1S_2U_1U_2}{Y_3},
$$
which follow from (\ref{Yi}), (\ref{1<2<3}) and (\ref{y3s3}), respectively.
But then it follows that
\begin{align*}
N^{(3)}(B)&\ll B \sum_{S_0,S_i,U_i,Y_1,Y_3}
\frac{S_1Y_1}{S_0S_2S_3U_2U_3Y_3}\\
&\ll B \sum_{S_0,S_i,U_i,Y_3}
\frac{S_3}{S_0S_1S_2U_1U_2}\\
&\ll B \sum_{S_0,S_1,S_2,U_i,Y_3} Y_3^{-1},
\end{align*}
whence $N^{(3)}(B)\ll B (\log B)^6$ if $Y_3U_3S_3^2 \gg
Y_2U_2S_2^2$.  Next we suppose that $Y_3U_3S_3^2 \ll
Y_2U_2S_2^2$, and take $\theta_3=1$ in (\ref{sf}).   
Observe that
$$
S_0\ll \frac{Y_3S_3}{S_1S_2U_1U_2}, \quad U_3 \ll
\frac{Y_2U_2S_2^2}{Y_3S_3^2},  \quad Y_1 \ll \frac{B}{Y_2Y_3},
$$
which follow from (\ref{beer}), the inequality $Y_3U_3S_3^2 \ll
Y_2U_2S_2^2$ and (\ref{Y}), respectively. We then argue as above to deduce that 
\begin{align*}
N^{(3)}(B)&=\sum_{S_0,S_i,U_i,Y_i}
S_0S_1S_2^{-1}S_3U_1U_3Y_1\\
&\ll\sum_{S_i,U_i,Y_i} \frac{Y_1Y_3U_3S_3^2}{U_2S_2^2}\\
&\ll\sum_{S_i,U_1,U_2,Y_i} Y_1Y_2\\
&\ll B \sum_{S_i,U_1,U_2,Y_2,Y_3}Y_3^{-1}.
\end{align*}
Hence $N^{(3)}(B)\ll B (\log B)^6$ in this case also.

Finally, we must estimate (\ref{sf}), whenever $i\neq 3$.  
Suppose that $i=1$, so that
$$
N^{(1)}(B)=\sum_{S_0,S_i,U_i,Y_i} \frac{S_0SU Y_2Y_3}{Y_3U_3S_3^2} =
\sum_{S_0,S_i,U_i,Y_i} S_0S_1S_2S_3^{-1}U_1U_2Y_2.
$$
Once again we separate our arguments according to the size of
$Y_3U_3S_3^2$.  Suppose that  $Y_3U_3S_3^2 \gg Y_2U_2S_2^2$, so that (\ref{y3s3})
holds.  Then $Y_3$ is fixed by the choices of $S_0,S_i,U_1,U_2$.
Moreover, we have the inequalities
$$
Y_2\ll \frac{B}{U_2US_0^2S_2^2}, \quad S_1 \ll \frac{S_0S_2S_3U_2U_3}{Y_1}, 
$$
which follow from (\ref{Yi}) and the fact that $Y_1U_1S_1^2 \ll 
Y_3U_3S_3^2 \ll S_0SU$, respectively.  Summing over $Y_2$, and then
over $S_1$, we therefore deduce that
\begin{align*}
N^{(1)}(B)
&\ll B \sum_{S_0,S_i,U_i,Y_1} \frac{S_1}{S_0S_2S_3U_2U_3}\\
&\ll B \sum_{S_0,S_2,S_3,U_i,Y_1} Y_1^{-1}\\
&\ll B (\log B)^6,
\end{align*}
in this case.

Alternatively, if $Y_3U_3S_3^2 \ll Y_2U_2S_2^2$, then $U_3$ is
determined by choices of $S_2,S_3,U_2,Y_2,Y_3$, and it follows that 
$$
N^{(1)}(B)\ll \sum_{S_0,S_i,U_1,U_2,Y_i} S_0S_1S_2S_3^{-1}U_1U_2Y_2.
$$
Upon summing over $S_0 \ll Y_3S_3/(S_1S_2U_1U_2)$, and then over $Y_2 \ll
B/(Y_{1}Y_3)$, we derive the estimate 
$$
N^{(1)}(B)\ll B\sum_{S_i,U_1,U_2,Y_1,Y_3} Y_1^{-1} \ll B (\log B)^6,
$$
in this case.   

An entirely similar argument handles the case $i=2$.  Upon
combining our various estimates we therefore deduce the statement of
Proposition \ref{n2}.

\end{document}